\documentclass[10pt]{amsart}

\usepackage{amsmath,amssymb,amsfonts,epsfig}

\newcommand{\rb}{\raisebox}
\newcommand{\ig}{\includegraphics}
\newcommand\risS[6]{\rb{#1pt}[#5pt][#6pt]{\begin{picture}(#4,15)(0,0)
  \put(0,0){\ig[width=#4pt]{#2.eps}} #3
     \end{picture}}}
\newcommand\smf[1]{\risS{-4}{#1}{}{15}{15}{8}}

\def\a{\alpha}
\def\b{\beta}
\def\d{\delta}
\def\ve{\varepsilon}
\def\Ga{\Gamma}
\def\cS{\mathcal S}
\def\cF{\mathcal F}
\def\bc{f}
\def\Z{{\mathbb Z}}
\newtheorem{thm}{Theorem}[section]
\newtheorem{defn}[thm]{Definition}
\newtheorem{exa}[thm]{Example}
\newtheorem{exas}[thm]{Examples}
\newtheorem{lemma}[thm]{Lemma}
\newtheorem{prop}[thm]{Proposition}
\newtheorem{cor}[thm]{Corollary}
\newtheorem{rem}[thm]{Remark}
\newtheorem{subs}[thm]{}   %\subsection

\def\kb#1{[ #1 ]}
\def\wo{\overline}
\def\ul{\underline}
\def\wt{\widetilde}
\def\plinc{\mbox{\begin{picture}(5,5)(0,0)
          \put(4,2){\circle{11}} \put(0,0){\mbox{$+$}}
                 \end{picture}}}
\def\miinc{\mbox{\begin{picture}(5,5)(0,0)
          \put(3,2){\circle{11}} \put(-.5,0){\mbox{$-$}}
                 \end{picture}}}
\def\vyd#1{\mbox{\underline{\it\bfseries #1}}\vspace{3pt}}
\def\pn#1#2#3{\put(#1,#2){\mbox{\tt\scriptsize #3}}}
\def\splinc{\mbox{\begin{picture}(5,5)(0,0)
          \put(4,2){\circle{7}} \put(1,.3){\mbox{\scriptsize $+$}}
                 \end{picture}}}
\def\smiinc{\mbox{\begin{picture}(5,5)(0,0)
          \put(3,2){\circle{7}} \put(-.2,.1){\mbox{\scriptsize $-$}}
                 \end{picture}}}
\begin{document}%%%%%%%%%%%%%%%%%%%%%%%%%%%%%%%%%%%%%%%%%%%%%%%%%%%%%%%%
%%%%%%%%%%%%%%%%%%%%%%%%%%%%%%%%%%%%%%%%%%%%%%%%%%%%%%%%%%%%%%%%%%%%%%%%

\title[Duality for graphs on surfaces and the Bollob\'as-Riordan polynomial]
  {Generalized duality for graphs on surfaces and the signed Bollob\'as-Riordan  
   polynomial}
\author[Sergei~Chmutov]{Sergei~Chmutov}
\date{}

\keywords{Graphs on surfaces, ribbon graphs, Bollob\'as-Riordan polynomial, Tutte polynomial, duality, virtual links, Jones polynomial.}

\begin{abstract}
We generalize the natural duality of graphs embedded into a surface to a duality with respect to a subset of edges. The dual graph might be embedded into a different surface. We prove a relation between the signed Bollob\'as-Riordan polynomials
of dual graphs. This relation unifies various recent results expressing the Jones polynomial of links as specializations of the Bollob\'as-Riordan polynomials.
\end{abstract}

\maketitle

\section*{Introduction} \label{s:intro}
We suggest a far reaching generalization of the famous duality relation,
$T_{\Ga}(x,y)=T_{\Ga*}(y,x)$, between the Tutte polynomials of a plane graph $\Ga$ and its natural dual graph $\Ga^*$ to graphs embedded into a higher genus surface. 

We generalize the duality to a duality with respect to a subset of edges. The dual graph might be embedded into a different surface. A flavor of our duality can be seen on the example (see details in example \ref{exdu}(b) on p.\pageref{exdu-b}).
$$\risS{0}{rg-torus}{\pn{15}{25}{$e$}
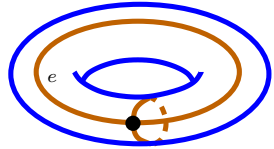
     \pn{20}{-10}{Graph $\Ga$ on a torus}}{100}{60}{25}\hspace{3cm}
\risS{5}{rg-sphere}{\pn{18}{30}{$e$}
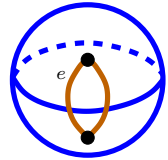
     \pn{-30}{-10}{Dual graph $\Ga^{\{e\}}$ with respect to the}
     \pn{-25}{-20}{edge $e$ is embedded into a sphere}}{60}{0}{0}\hspace{2cm}
$$

In this paper we are dealing with {\it cellularly} embedded graphs. Such embeddings can be formalized in a notion of {\it ribbon graphs}.
A {\it ribbon graph} $G$ is a surface with boundary and a decomposition into a union of closed topological discs of two types, vertices and edges, subject to some natural axioms (see the precise definition below).
If we shrink each vertex-disc to its center and each edge-disc to a line connecting the central points of its vertices, we will get an ordinary graph $\Ga$, the {\it underlying graph}, embedded into the surface of $G$. Conversely, any graph $\Ga$ embedded into a surface determines a ribbon graph structure on a small neighborhood $G$ of $\Ga$ inside the surface. Thus ribbon graphs are nothing else than abstract graphs cellularly embedded into a closed surface.

For any ribbon graph $G$, there is a {\it natural} dual ribbon graph $G^*$, also called
{\it Euler-Poincar\'e} dual, defined as follows. First we glue a disc, {\it face}, to each boundary component of $G$, obtaining a closed surface $\wt{G}$ without boundary. 
Then we remove the interior of all vertex-discs of $G$.
The newly glued discs-faces will be the vertex-discs of $G^*$. The edge-discs for $G^*$ will be the same as for $G$ but their attachments to new vertices will, of course, be different.
Both underlying graphs $\Ga$ and $\Ga^*$ turn out to be embedded into the same surface
$\wt{G}=\wt{G^*}$ in a natural dual manner: the edges of $\Ga$ are mutually perpendicular to the edges of $\Ga^*$, the vertices of $\Ga$ correspond to the faces of $\Ga^*$ and vise versa.
 
B.~Bollob\'as and O.~Riordan \cite{BR3} found a generalization of the Tutte polynomial
for ribbon graphs which captures some topological information. For non-planar graphs, there is no duality relation for the Tutte polynomial but there is one for the Bollob\'as-Riordan polynomial. In \cite{BR3}, it was proved for one free variable. J.~Ellis-Monaghan and I.~Sarmiento \cite{EMS} extended it to a two free variables relation (see also \cite{Mof}). 

We work with {\it signed ribbon graphs}, that is ribbon graphs whose edges are marked by either $+1$ or $-1$. For such graphs, we generalize the notion
of duality to the {\it duality with respect to a subset of edges}. Let $E'\subseteq E(G)$
be a subset of edges of a ribbon graph $G$. The dual graph $G^{E'}$ is constructed in the following way. Consider the spanning subgraph $F_{E'}$ of $G$ containing all the vertices of $G$ and the edges from $E'$ only. Glue a disc-face into each boundary component of $F_{E'}$;
these faces of $F_{E'}$ are going to be the vertex-discs of the dual graph $G^{E'}$. Removing the interior of all old vertices of $G$ we get $G^{E'}$. Its edges are the same discs as in $G$ only the attachments of edges from $E'$ to new vertices are changed. The signs of edges in $E'$ have to be changed to the opposite. In general, the genus of $G^{E'}$ is not equal to the genus of $G$. So the corresponding underlying graphs are embedded into different surfaces.

We give a duality relation for the signed version of the Bollob\'as-Riordan polynomial (introduced in \cite{CP}) of graphs $G$ and $G^{E'}$. When $E'=E(G)$ and all edges of $G$ are positive, our relation essentially coincides with the one from \cite{EMS,Mof}. If, moreover, $G$ is planar then our duality relation reduces to the famous duality relation for the Tutte polynomial. 

Igor Pak suggested to use the Bollob\'as-Riordan polynomial in knot theory for  Thistlethwaite's type theorems. This was first realized in \cite{CP}. Then there were
two other realizations of this idea in \cite{DFKLS,CV}. Formally all three theorems from
\cite{CP,DFKLS,CV} are different. They used different constructions of a ribbon graph from a link diagram and different substitutions in the Bollob\'as-Riordan polynomials of these graphs. Here we show that our duality relation allows to derive the theorems of \cite{DFKLS,CP} from the one of \cite{CV}. 

The generalized duality and the corresponding relation for the Bollob\'as-Riordan polynomial appeared in my attempts to unify the theorems from \cite{CP,DFKLS,CV}. 
I express my deep gratitude to their authors. Also I am thankful to M.~Chmutov for
useful discussions and for showing me the signed version of the contraction-deletion property (proposition \ref{prop:c-d}) for the Bollob\'as-Riordan polynomial and to
I.~Moffatt and F. Vignes-Tourneret for very interesting comments.

\bigskip
\section{Ribbon graphs and generalized duality}\label{s:rg}

As said above, ribbon graphs are practically the same as graphs on surfaces and thus
they are objects of Topological Graph Theory \cite{GT,LZ,MT}. From the point of view of this theory, our ribbon graphs are nothing else but
the band decompositions from \cite[section 3.2]{GT} with the interiors of all 2-bands removed. Oriented ribbon graphs appear under different names such as {\it rotation systems} \cite{MT}, {\it maps}, {\it fat graphs}, {\it cyclic graphs}, {\it dessins d'Enfants} \cite{LZ}. Since the pioneering paper of L.~Heffter of 1891 \cite{He} they occur in various parts of mathematics ranging from graph theory, combinatorics, and topology to representation theory, Galois theory, algebraic geometry, and quantum field theory \cite{CDBooK,Ha,Ko,LZ,Pe,RT}. 

\noindent\parbox[t]{4.5in}{For example, ribbon graphs are used to enumerate cells in the cell decomposition of the moduli spaces of complex algebraic curves
\cite{Ha,Pe,LZ}. The absolute Galois group $\mathrm{Aut}(\wo{\mathbb Q}/{\mathbb Q})$ faithfully acts on the set of ribbon graphs (see \cite{LZ} and references therein). Ribbon graphs are the main combinatorial objects of the Vassiliev knot invariant theory
\cite{CDBooK}. They are very useful for Hamiltonicity of the Cayley graphs \cite{GM,GY}. The ribbon graph on the right represents the Cayley graph of the
$\langle 5,3,2\rangle$ presentation\linebreak }\qquad
\parbox[t]{1.5in}{$\risS{-65}{grap-1}{}{70}{0}{0}$}\vspace{-8pt}\\
of the alternating group $A_5 = \langle x,y\ \bigr|\ x^5=y^3=(xy)^2=1\rangle$ (that is also isomorphic to $\mathrm{PSL}_2(5)$), where $x=(1\, 2\, 3\, 4\, 5)$ and 
$y=(2\, 5\, 4)$, see \cite[Fig.1]{GY}.

\newpage
We will use a formal definition from \cite{BR3}. 

\begin{defn}\label{def:rb}{\rm
A {\it ribbon graph} $G$ is a surface (possibly non-orientable) with boundary, represented as the union of two sets of closed topological discs called
{\it vertices} $V(G)$ and {\it edges} $E(G)$, satisfying the following
conditions:%\vspace{-2pt}
\begin{itemize}
\item[$\bullet$] these vertices and edges intersect by disjoint line segments;
\item[$\bullet$] each such line segment lies on the boundary of precisely one vertex and precisely one edge;
\item[$\bullet$] every edge contains exactly two such line segments.
\end{itemize} 

\noindent
A ribbon graph is said to be {\it signed} if it is accompanied with a {\it sign function} $\ve: E(G) \to \{\pm 1\}$.}
\end{defn}

We consider ribbon graphs up to a homeomorphism of the corresponding surfaces preserving the decomposition on vertices and edges.

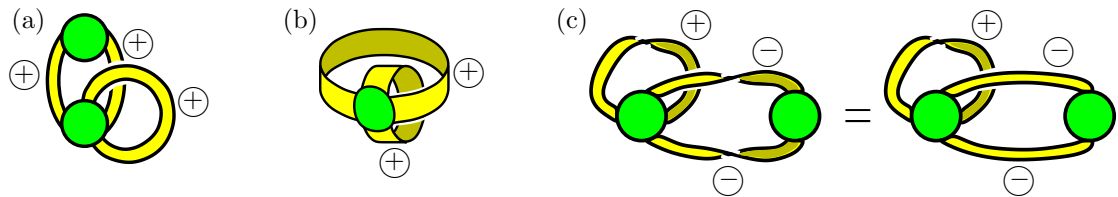
\begin{exas}\label{ex}{\rm %(If the sign is omitted it is assumed to be $+1$.)
$$\mbox{(a)}\risS{-52}{rg-ex1}{\put(-12,32){\plinc} \put(31,44){\plinc} 
           \put(52,22){\plinc}}{50}{0}{0}\hspace{1.4cm}
  \mbox{(b)}\risS{-44}{rg-ex2}{\put(53,25){\plinc} \put(25,-9){\plinc}
           }{50}{0}{0}\hspace{1.4cm}
  \mbox{(c)}\risS{-52}{rg-ex3}{\put(37,51){\plinc}\put(65,41){\miinc}\put(50,-8){\miinc}
      \put(147,51){\plinc}\put(175,41){\miinc}\put(160,-8){\miinc}}{200}{0}{65}
$$}
\end{exas}

It is important to note that a ribbon graph is an abstract two-dimensional surface with boundary; its embedding into the 3-space is irrelevant.

\bigskip
It may be convenient to have a more combinatorial definition of ribbon graphs.
We may think about an edge not as a disc, but rather as a rectangle attached to the corresponding vertices along a pair of its opposite sides. Pick an orientation for each vertex-disc and for each edge-rectangle and label the edges. The orientations of the rectangles induce arrows on their sides. Then we draw all vertex-discs as disjoint circles in the plane oriented counterclockwise,
but instead of drawing edges we draw only the arrows of the corresponding sides on the boundary circles of vertices and put the corresponding labels. Here is an illustration of this procedure for the example \ref{ex}(a). 
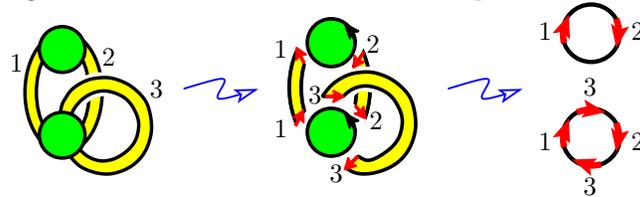
$$\risS{-48}{rg-ex1}{\put(-5,40){\mbox{$1$}} 
         \put(30,42){\mbox{$2$}}\put(48,30){\mbox{$3$}}}{50}{0}{0}\quad
  \risS{-20}{toto}{}{30}{0}{0}\quad
  \risS{-48}{rgs-exa1}{\put(-5,45){\mbox{$1$}}\put(-5,15){\mbox{$1$}} 
         \put(30,47){\mbox{$2$}}\put(31,18){\mbox{$2$}} 
         \put(16,-2){\mbox{$3$}}\put(8,28){\mbox{$3$}}}{50}{0}{0}\quad
  \risS{-20}{toto}{}{30}{0}{0}\quad
  \risS{-48}{rgs-exa2}{\put(-5,52){\mbox{$1$}}\put(-5,10){\mbox{$1$}} 
         \put(30,52){\mbox{$2$}}\put(30,10){\mbox{$2$}} 
         \put(12,-7){\mbox{$3$}}\put(12,31){\mbox{$3$}}}{30}{10}{55}\quad
$$
The resulting figure uniquely determines the ribbon graph. In the case of signed ribbon graphs, besides circles with arrows, we need a list of edge labels with the signs. Here are the figures corresponding to examples \ref{ex}(b,c).
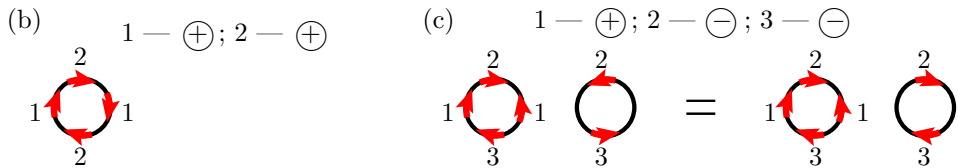
$$\mbox{(b)}\risS{-45}{rgs-exb1}{\put(-5,10){\mbox{$1$}}\put(30,10){\mbox{$1$}} 
         \put(12,-7){\mbox{$2$}}\put(12,31){\mbox{$2$}}
      \put(30,40){\mbox{1 ---\ \ \plinc\ \ ; 2 ---\ \ \plinc}}}{30}{5}{60} \hspace{4cm}
\mbox{(c)}\risS{-45}{rgs-exc1}{\put(-5,10){\mbox{$1$}}\put(30,10){\mbox{$1$}}
         \put(12,30){\mbox{$2$}}\put(53,30){\mbox{$2$}} 
         \put(12,-7){\mbox{$3$}}\put(53,-7){\mbox{$3$}}
         \put(117,10){\mbox{$1$}}\put(152,10){\mbox{$1$}}
         \put(134,30){\mbox{$2$}}\put(175,30){\mbox{$2$}} 
         \put(134,-7){\mbox{$3$}}\put(175,-7){\mbox{$3$}}
         \put(30,45){\mbox{1 ---\ \ \plinc\ \ ; 2 ---\ \ \miinc\ \ ; 3 ---\ \ \miinc}}
         }{190}{5}{60}\quad
$$
The last example shows that we can reverse all arrows on a circle together with reversing the cyclic order of the arrows along the circle. This would correspond to picking the opposite orientation of the vertex-disc. Also we can reverse two arrows with the same label. This would correspond to picking the opposite orientation of the edge-rectangle. %Now we can define ribbon graphs as follows.

\begin{defn} {\bf Arrow presentation.\ }\label{def:rb2}{\rm
A {\it ribbon graph} $G$ is a collection of disjoint (non-nested) oriented circles in the plane with a bunch of
labeled arrows on them, such that each label occurs precisely twice. Two such collections are considered equivalent if one is obtained from the other by reversing all arrows on a circle and reversing the cyclic order of the arrows along it, or by reversing two arrows with the same labels.
A ribbon graph is said to be {\it signed} if it is accompanied by a {\it sign function} from the set of labels to $\{\pm 1\}$.}
\end{defn}

The information from this definition allows to recover the ribbon graph in the sense of definition \ref{def:rb}. To construct the vertices we have to fill in the circles; the labels indicate the edges which have to be glued to the vertices;
and an arrow determines the places of gluing and the character of gluing, that is whether we should glue an edge as a planar band or we should make a twist on it (the head of an arrow should go along a side of the edge and come to the tail of another arrow with the same label).

%\bigskip
\begin{defn}\label{def:gedu}{\rm
Let $E'\subseteq E(G)$ be a subset of edges of a signed ribbon graph $G$. We define the {\it dual graph $G^{E'}$
with respect to the subset $E'$} as follows. Consider the spanning subgraph $F_{E'}$ of $G$ containing all the vertices of $G$ and only the edges from $E'$. Each boundary component of $F_{E'}$ will be a circle from the collection for $G^{E'}$. The arrows from the edges of $E(G)\setminus E'$ will be the same as for $G$ itself. The edges
from $E'$ will give different arrows. Let $e\in E'$. The rectangle representing $e$ intersects with vertex-discs of $G$ by a pair of opposite sides. But it
intersects with the boundary of the surface $F_{E'}$ by another pair of its opposite sides. This second pair gives a pair of arrows on the circles of $G^{E'}$ corresponding to $e$.

\noindent
The sign function $\ve_{G^{E'}}$ for $G^{E'}$ is defined by the equations: 
$\ve_{G^{E'}}(e)=-\ve_G(e)$ if $e\in E'$, and 
$\ve_{G^{E'}}(e)=\ve_G(e)$ if $e\not\in E'$.}
\end{defn}

%\bigskip
\begin{subs}\label{ss:dus} {\bf Duality from the surface point of view.}\end{subs}

Suppose an edge $e$ connects two different vertices and it is attached to them along its sides $a$ and $c$. To construct the dual graph $G^{\{e\}}$ we double this edge so that the two copies will be sewed together along the sides $b$ and $d$ forming a cylinder.
Then we regard the union of the two vertices and one of the copies of $e$ as a new vertex. It will be the vertex of $G^{\{e\}}$. Now the second copy of $e$ turns out to be attached to this new vertex along the sides $b$ and $d$. It constitutes the edge $e$
of the dual graph $G^{\{e\}}$. 
So the dual graph $G^{\{e\}}$ may be considered as a contraction $G/e$ with an additional edge, the second copy of $e$.
Pictorially this looks as follows (the boxes with dashed arcs mean that there might be other edges attached to these vertices).
$$G =\ 
  \risS{-8}{dus1}{\put(27,10){\mbox{$a$}}\put(40,10){\mbox{$c$}}
          \put(34,-4){\mbox{$b$}}\put(34,21){\mbox{$d$}}}{70}{0}{0}\quad 
   \risS{-4}{toto}{\put(5,15){\mbox{\scriptsize\tt double}}
         \put(0,-8){\mbox{\scriptsize\tt the edge}}}{30}{0}{0}\quad 
  \risS{-35}{dus2}{}{110}{0}{0}\quad 
   \risS{-4}{toto}{\put(-5,15){\mbox{\scriptsize\tt forming the}}
         \put(-5,-8){\mbox{\scriptsize\tt new vertex}}}{30}{0}{0}\quad 
  \risS{-35}{dus3}{\put(70,23){\mbox{$b$}}\put(70,48){\mbox{$d$}}}{110}{45}{30} 
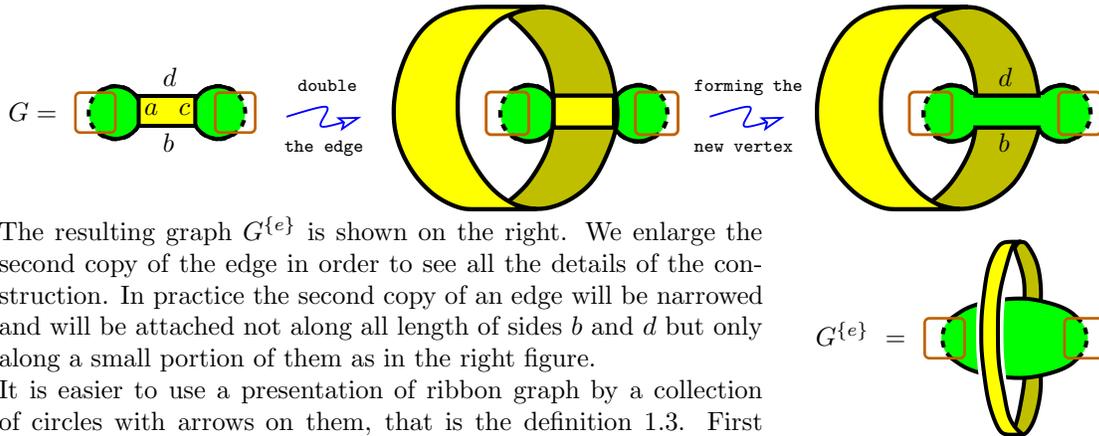
$$
\parbox[t]{4in}{The resulting graph $G^{\{e\}}$ is shown on the right. We enlarge the second copy of the edge in order to see all the details of the construction.
In practice the second copy of an edge will be narrowed and will be attached not along all length of sides $b$ and $d$ but only along a small portion of them as in the right figure. 

It is easier to use a presentation of ribbon graph by a collection of circles with arrows on them, that is the definition \ref{def:rb2}. First we pick  orientations on edges of $E'$. In our current~example,\linebreak }\qquad
\parbox[t]{2in}{\rb{-40pt}{$G^{\{e\}}\ =\ \risS{-35}{dus4}{}{70}{10}{0}$}}
              \vspace{-12pt}\\
$E'=\{e\}$, and we may pick the orientation $a$ --- $b$ --- $c$ --- $d$ and indicate it by the arrows. Then we draw the boundary components of $F_{E'}$, keeping on them only the arrows of sides which were free, not attached to the vertices of $G$. In our case, these are the sides $b$ and $d$. After that, we straighten the boundary components into circles and pull them apart. What we obtain will be a presentation of the graph $G^{\{e\}}$. Here is how it works.
$$G =\ 
  \risS{-8}{dus5}{}{70}{0}{0}\hspace{30pt}
   \risS{-4}{toto}{\put(-5,15){\mbox{\scriptsize\tt drawing the}}
     \put(-15,-8){\mbox{\scriptsize\tt boundary components}}}{30}{0}{0}\hspace{30pt}
  \risS{-8}{dus6}{}{70}{0}{0}\quad=\quad
  \risS{-30}{dus7}{\put(30,10){\mbox{$e$}}\put(30,52){\mbox{$e$}}}{70}{20}{30} 
  \quad=\ G^{\{e\}}
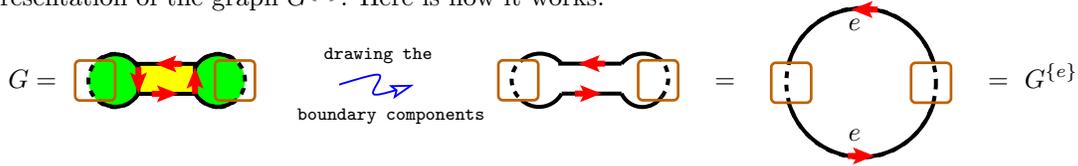
$$

\begin{exas}\label{exdu}{\rm
These are the continuations of examples \ref{ex}. They demonstrate the construction of dual graphs with respect to some subsets of edges.
$$\rb{80pt}{(a)}\qquad\begin{array}{l}
G =\ \risS{-25}{rg-ex1}{\put(40,55){\mbox{$E'=\{3\}$}}\put(-5,40){\mbox{$1$}} 
         \put(30,42){\mbox{$2$}}\put(48,30){\mbox{$3$}}
     \put(-40,70){\mbox{1 ---\ \ \plinc\ \ ; 2 ---\ \ \plinc\ \ ; 3 ---\ \ \plinc}}
         }{50}{60}{70}\qquad 
   \risS{0}{toto}{}{30}{15}{8}\qquad 
   \risS{-25}{du-ex1}{\put(-5,45){\mbox{$1$}}\put(-5,15){\mbox{$1$}} 
         \put(25,45){\mbox{$2$}}\put(25,12){\mbox{$2$}} 
         \put(32,20){\mbox{$3$}}\put(50,15){\mbox{$3$}} 
         \put(0,65){\mbox{$\partial F_{\{3\}}$}}
         \put(50,70){\mbox{1 ---\ \ \plinc\ \ ; 2 ---\ \ \plinc\ \ ; 3 ---\ \ \miinc}}
         }{50}{0}{0}\qquad = \qquad 
   \risS{-25}{du-ex1a}{\put(-5,42){\mbox{$1$}}\put(17,7){\mbox{$1$}} 
         \put(27,44){\mbox{$2$}}\put(48,7){\mbox{$3$}} 
         \put(38,44){\mbox{$2$}}\put(70,45){\mbox{$3$}}}{70}{0}{0}\quad= \\ \hspace{3.5cm}
=\qquad\risS{-25}{du-ex1b}{\put(-5,42){\mbox{$1$}}\put(17,7){\mbox{$1$}} 
         \put(27,44){\mbox{$2$}}\put(48,7){\mbox{$3$}} 
         \put(38,44){\mbox{$2$}}\put(70,45){\mbox{$3$}}}{70}{0}{0}\quad=\qquad 
   \risS{-27}{du-ex1c}{\put(1,26){\plinc}
         \put(30,56){\plinc}\put(55,26){\miinc}}{65}{0}{50}\quad=\quad G^{\{3\}}
\end{array}
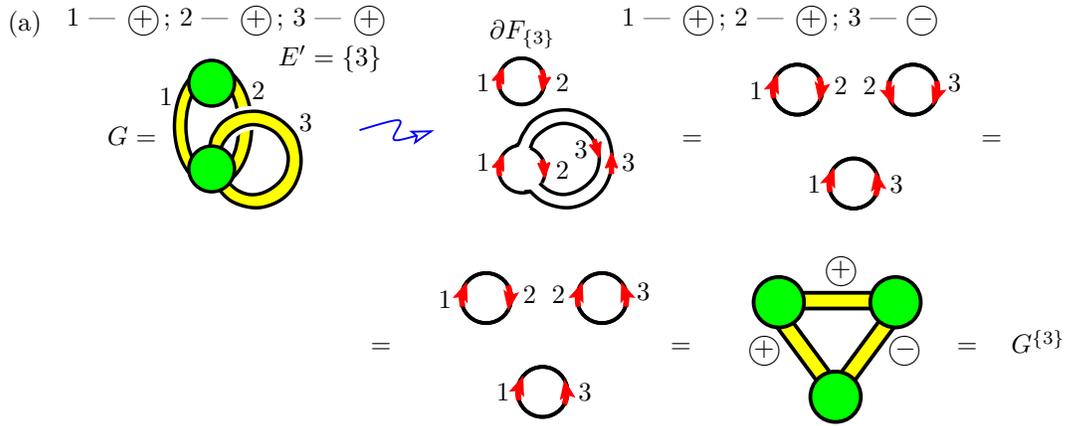
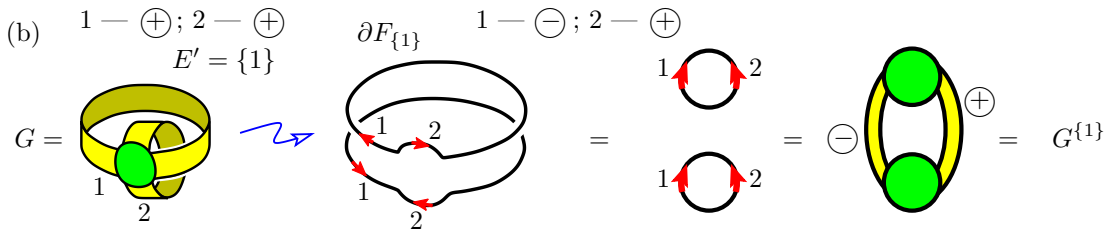
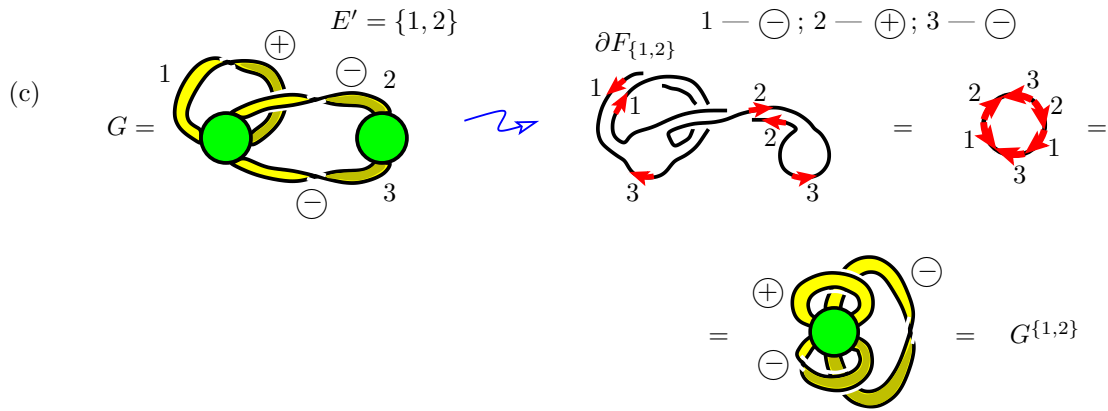
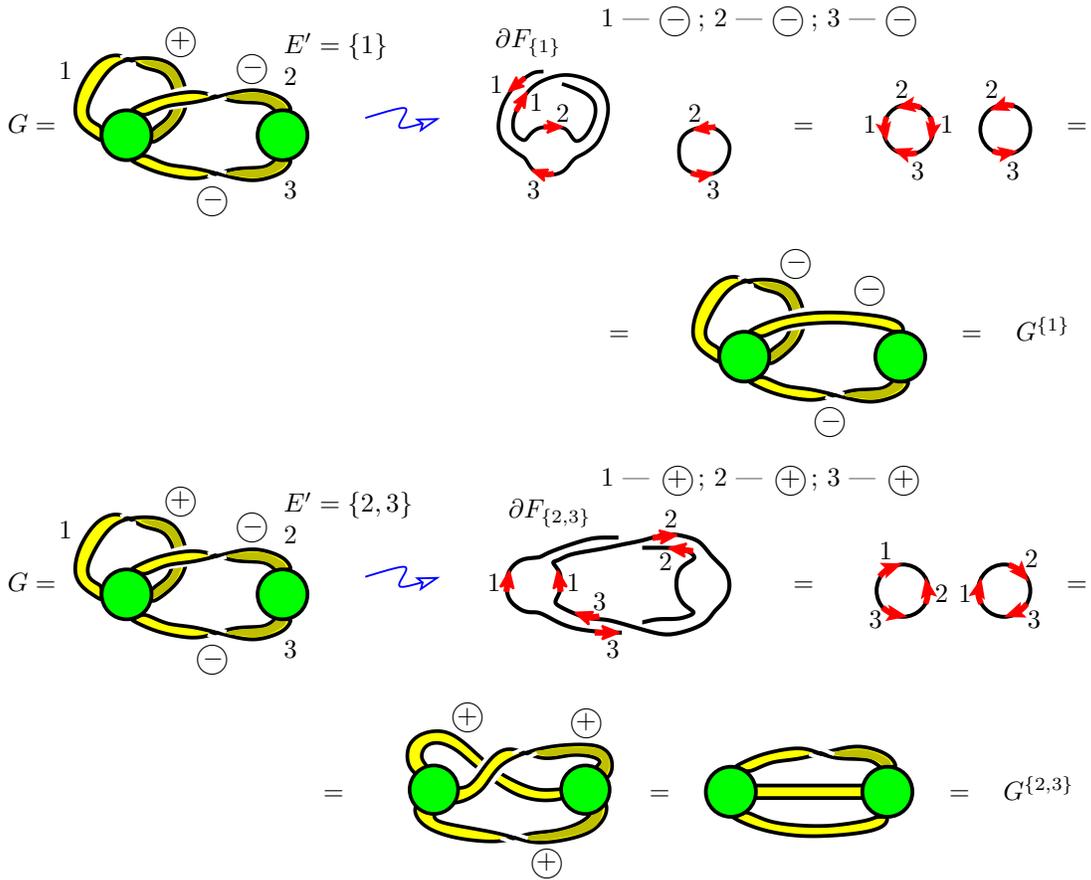
$$
$$\rb{40pt}{(b)}\hspace{-10pt}\label{exdu-b}
G =\ \risS{-20}{rg-ex2}{\put(35,50){\mbox{$E'=\{1\}$}}
         \put(5,2){\mbox{$1$}}\put(22,-8){\mbox{$2$}}
         \put(0,65){\mbox{1 ---\ \ \plinc\ \ ; 2 ---\ \ \plinc}} }{50}{10}{70}\quad 
   \risS{0}{toto}{}{30}{0}{0}\quad 
   \risS{-25}{du-ex2}{\put(12,30){\mbox{$1$}}\put(5,5){\mbox{$1$}} 
         \put(32,28){\mbox{$2$}}\put(25,-6){\mbox{$2$}} 
         \put(5,65){\mbox{$\partial F_{\{1\}}$}}
         \put(50,70){\mbox{1 ---\ \ \miinc\ \ ; 2 ---\ \ \plinc}}
         }{70}{0}{0}\qquad = \qquad
  \risS{-25}{du-ex2a}{\put(-5,52){\mbox{$1$}}\put(-5,10){\mbox{$1$}} 
         \put(30,52){\mbox{$2$}}\put(30,10){\mbox{$2$}}}{30}{0}{0}\quad=\qquad 
   \risS{-25}{du-ex2b}{\put(-11,26){\miinc}\put(40,40){\plinc}}{38}{0}{00}
\quad=\quad G^{\{1\}}
$$
$$\rb{50pt}{(c)}\qquad\begin{array}{l}\label{p:ex-c}
G =\ \risS{-20}{rg-ex31}{\put(60,60){\mbox{$E'=\{1,2\}$}}
         \put(-5,40){\mbox{$1$}}\put(80,38){\mbox{$2$}}\put(80,-5){\mbox{$3$}}
         \put(37,52){\plinc}\put(65,42){\miinc}\put(50,-8){\miinc}}{90}{50}{70}\qquad 
   \risS{0}{toto}{}{30}{15}{8}\qquad 
\risS{-20}{du-ex31}{\put(-2,35){\mbox{$1$}}\put(13,28){\mbox{$1$}}
         \put(60,33){\mbox{$2$}}\put(64,16){\mbox{$2$}} 
         \put(12,-5){\mbox{$3$}}\put(80,-5){\mbox{$3$}}
         \put(0,52){\mbox{$\partial F_{\{1,2\}}$}}
         \put(40,60){\mbox{1 ---\ \ \miinc\ \ ; 2 ---\ \ \plinc\ \ ; 3 ---\ \ \miinc}}
         }{90}{0}{0}\qquad=\qquad 
  \risS{-12}{du-ex31a}{\put(15,-6){\mbox{$3$}}\put(-4,7){\mbox{$1$}} 
         \put(-2,25){\mbox{$2$}}\put(20,30){\mbox{$3$}}
         \put(29,20){\mbox{$2$}}\put(28,4){\mbox{$1$}}}{30}{0}{0}\quad=\qquad \\
\hspace{8cm}=\qquad 
   \risS{-25}{du-ex31b}{\put(-9,15){\miinc}\put(-12,42){\plinc}
         \put(49,50){\miinc}}{50}{0}{40}
\quad=\quad G^{\{1,2\}}
\end{array}
$$
$$\rb{50pt}{\ }\begin{array}{l}
G =\ \risS{-20}{rg-ex31}{\put(80,50){\mbox{$E'=\{1\}$}}
\put(-5,40){\mbox{$1$}}\put(80,38){\mbox{$2$}}\put(80,-5){\mbox{$3$}}
         \put(37,52){\plinc}\put(65,42){\miinc}\put(50,-8){\miinc}}{90}{40}{70}\qquad 
   \risS{0}{toto}{}{30}{65}{8}\qquad 
\risS{-20}{du-ex32}{\put(-2,35){\mbox{$1$}}\put(13,28){\mbox{$1$}}
         \put(23,24){\mbox{$2$}}\put(73,24){\mbox{$2$}} 
         \put(12,-5){\mbox{$3$}}\put(80,-5){\mbox{$3$}}
         \put(0,52){\mbox{$\partial F_{\{1\}}$}}
         \put(40,60){\mbox{1 ---\ \ \miinc\ \ ; 2 ---\ \ \miinc\ \ ; 3 ---\ \ \miinc}}
         }{90}{0}{0}\qquad=\qquad 
  \risS{-12}{du-ex32a}{\put(14,-6){\mbox{$3$}}\put(-4,12){\mbox{$1$}} 
         \put(8,25){\mbox{$2$}}\put(50,-6){\mbox{$3$}}
         \put(42,25){\mbox{$2$}}\put(25,12){\mbox{$1$}}}{60}{0}{0}\quad=\qquad \\
\hspace{8cm}=\qquad 
   \risS{-25}{du-ex32b}{\put(37,52){\miinc}\put(65,42){\miinc}
              \put(50,-8){\miinc}}{90}{0}{50}\quad=\quad G^{\{1\}}
\end{array}
$$
$$\rb{50pt}{\ }\begin{array}{l}
G =\ \risS{-20}{rg-ex31}{\put(80,50){\mbox{$E'=\{2,3\}$}}
\put(-5,40){\mbox{$1$}}\put(80,38){\mbox{$2$}}\put(80,-5){\mbox{$3$}}
         \put(37,52){\plinc}\put(65,42){\miinc}\put(50,-8){\miinc}}{90}{40}{70}\qquad 
   \risS{0}{toto}{}{30}{15}{8}\qquad 
\risS{-20}{du-ex33}{\put(-3,20){\mbox{$1$}}\put(27,20){\mbox{$1$}}
         \put(62,28){\mbox{$2$}}\put(64,44){\mbox{$2$}} 
         \put(42,-5){\mbox{$3$}}\put(37,13){\mbox{$3$}}
         \put(5,48){\mbox{$\partial F_{\{2,3\}}$}}
         \put(40,60){\mbox{1 ---\ \ \plinc\ \ ; 2 ---\ \ \plinc\ \ ; 3 ---\ \ \plinc}}
         }{90}{0}{0}\qquad=\qquad 
  \risS{-12}{du-ex33a}{\put(-2,-2){\mbox{$3$}}\put(2,24){\mbox{$1$}} 
         \put(23,8){\mbox{$2$}}\put(58,-2){\mbox{$3$}}
         \put(57,21){\mbox{$2$}}\put(32,8){\mbox{$1$}}}{60}{0}{0}\quad=\qquad \\
\hspace{4.2cm}=\qquad 
   \risS{-18}{du-ex33b}{\put(20,47){\plinc}\put(65,45){\plinc}
              \put(50,-8){\plinc}}{80}{0}{0}\quad=\quad 
   \risS{-15}{du-ex33c}{}{80}{0}{40}\quad=\quad G^{\{2,3\}}
\end{array}
$$
}\end{exas}

%\newpage
\begin{subs}\label{ss:due} {\bf Duality with respect to a single edge.}\end{subs}
An important special case of duality is $E'=\{e\}$. We are going to figure out how to represent
the duality in terms of the collection or circles with arrows, i.e. in terms of definition
\ref{def:rb2}. There are three cases:
\begin{itemize}
\item[(i)] the edge $e$ connects two different vertices;
\item[(ii)] the edge $e$ connects a vertex with itself in an orientable way;
\item[(iii)] the edge $e$ connects a vertex with itself in a non-orientable way.
\end{itemize}

\ul{\bf Case (i).} 
We may choose the orientation on the vertices so that both arrows corresponding to the edge $e$ will point counterclockwise. It is easy to see that the duality will look as follows.
$$G\ =\ 
  \risS{-8}{due-i1}{\put(31,8){\mbox{$e$}}\put(37,8){\mbox{$e$}}}{70}{20}{0}\quad 
   \risS{-4}{toto}{}{30}{15}{8}\quad 
  \risS{-17}{due-i2}{\put(31,8){\mbox{$e$}}\put(37,27){\mbox{$e$}}}{70}{0}{15}\ 
  =\ G^{\{e\}}\ ,
$$
where the boxes $A$ and $B$ with dashed arcs mean that there might be other arrows attached to the vertices. 
Also there might be other circles with arrows on them, but they are the same
for both figures and therefore are omitted.

\medskip
\ul{\bf Case (ii).} This is precisely the opposite case to the previous one.
$$G\ =\ 
  \risS{-17}{due-i2}{\put(31,8){\mbox{$e$}}\put(37,27){\mbox{$e$}}}{70}{0}{15}\quad 
   \risS{-4}{toto}{}{30}{15}{8}\quad 
  \risS{-8}{due-i1}{\put(31,8){\mbox{$e$}}\put(37,8){\mbox{$e$}}}{70}{20}{0}\ 
  =\ G^{\{e\}}\ .
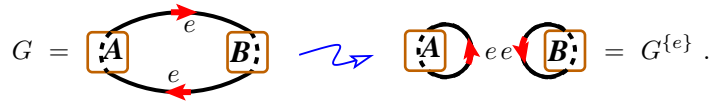
$$

\medskip
\ul{\bf Case (iii).} In this case the arrows $e$ point in opposite directions in the original graph $G$.
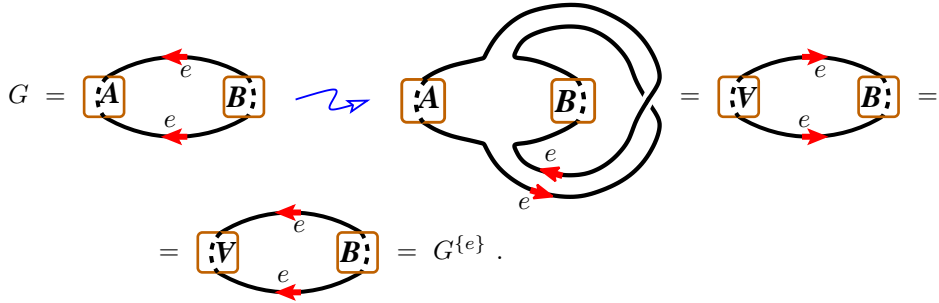
$$\begin{array}{l} G\ =\ 
  \risS{-17}{due-iii1}{\put(31,8){\mbox{$e$}}\put(37,27){\mbox{$e$}}}{70}{0}{0}\quad 
   \risS{-4}{toto}{}{30}{15}{8}\quad 
  \risS{-38}{due-iii2}{\put(45,-2){\mbox{$e$}}\put(55,16){\mbox{$e$}}}{100}{40}{50}\ 
  =\ \risS{-17}{due-iii3}{\put(31,8){\mbox{$e$}}\put(37,27){\mbox{$e$}}}{70}{0}{0}\ =\\
\hspace{2cm}
  =\ \risS{-17}{due-iii4}{\put(31,8){\mbox{$e$}}\put(37,27){\mbox{$e$}}}{70}{0}{20}\  
  =\ G^{\{e\}}\ .
\end{array}
$$
Here the upside-down box $A$ means that we should cut an arc with arrows inside the box $A$ of the original graph, flip it, and sew it back. The resulting graph will be $G^{\{e\}}$.

\bigskip
\begin{subs}\label{le:du-prop} {\bf Simple properties of duality.}\end{subs}
The next lemma describes some properties of the generalized duality.

{\bf Lemma.}{\it 
\begin{itemize}
\item[(a)] Suppose that an edge $e$ does not belong to $E'$. Then 
$G^{E'\cup\{e\}} = \bigl(G^{E'}\bigr)^{\{e\}}$.
\item[(b)] $\bigl(G^{E'}\bigr)^{E'}= G$.
\item[(c)] $\bigl(G^{E'}\bigr)^{E''}= G^{\Delta(E',E'')}$, where 
    $\Delta(E',E''):=(E'\cup E'')\setminus (E'\cap E'')$ is the symmetric difference of  
    sets.
\item[(d)] The generalized duality preserves orientability of ribbon graphs.
\item[(e)] Let $\wt{G}$ be a surface without boundary obtained from $G$ by gluing discs to all boundary component of $G$. Then $\wt{G^{E'}} = \wt{G^{E(G)\setminus E'}}$.
\item[(f)] The generalized duality preserves the number of connected components of
ribbon graphs.
\end{itemize}}

{\bf Proof.}

(a). This is a direct consequence of definition \ref{def:gedu}. It allows to find the dual with respect to $E'$ one edge at a time.

(b). By (a), we may assume that $E'$ consists of one edge, $E'=\{e\}$. The cases (i) --- (iii) above clearly show that the duality with respect to an edge is an involution.

(c). This follows from (a) and (b).

(d). In terms of the definition \ref{def:rb2}, orientability of the surface of a ribbon graph means that we can choose all the arrows coherently with the counterclockwise orientation of the circles. The figures of the cases (i) and (ii) above show that
if $G$ was orientable than $G^{\{e\}}$ will be orientable too. Now the statement follows from (a) and (b).

(e). Because of (a) and (b) it is sufficient to prove $\wt{G} = \wt{G^{E(G)}}$. This is obvious because the vertices of $G^{E(G)}$ are precisely the discs glued to $G$ to obtain 
$\wt{G}$ and vise versa. In particular, for a planar graph $G$, $G^{E(G)}=G^*$ is the dual planar graph in the ordinary sense.

(f). This is obvious.
\hspace*{\fill}$\square$

\bigskip
It is a consequence of the lemma that duality can be understood as an action of the group $\Z_2^l$ on ribbon graphs with $l$ edges. Then the number of non-isomorphic graphs (ignoring the sign function) dual to a given ribbon graph $G$ can be regarded an invariant of $G$. For the ribbon graph of example \ref{ex}(a) it is equal to 4;
for \ref{ex}(b) --- 2; for \ref{ex}(c) --- 5. Some of these duals are shown in example
\ref{exdu}. We leave it up to the reader to find the remaining duals as a useful exercise.

\newpage
%\bigskip
\begin{subs}\label{ss:du-c/d} {\bf Duality and contraction-deletion.}\end{subs}
For planar graphs (and more generally for {\it matroids}) it is well known that a contraction of an edge corresponds to a deletion of the edge in the dual graph. We can extend this property to ribbon graphs. Moreover, for ribbon graphs we can give a more subtle definition of a contraction of a (not necessary orientable) loop.

A \vyd{deletion} $G-e$ of an edge $e$ from a ribbon graph $G$ can be defined naturally as a graph obtained from $G$ by removal of the corresponding ribbon. In the arrow presentation (definition \ref{def:rb2}) it is reduced to the deletion of the pair of arrows labeled by $e$. 

A \vyd{contraction} $G/e$ is defined by the equation \fbox{$G/e\ :=\ G^{\{e\}}-e$}\ .

It is useful to consider the three cases from section \ref{ss:due}.

\noindent\parbox[t]{2in}{\ul{\bf Case (i).}\vspace{3pt} {\it The edge $e$ connects two different vertices.} \vspace{8pt}}\qquad
\parbox{4in}{$G = 
  \risS{-8}{due-i1}{\put(31,8){\mbox{$e$}}\put(37,8){\mbox{$e$}}}{70}{0}{0}\quad 
   \risS{-4}{toto}{}{30}{0}{0}\quad 
  \risS{-13}{contr-i2}{}{70}{25}{10} = G/e\ .$}
This is the familiar contraction of a non-loop.\vspace{-15pt} It is also called
{\it Whitehead collapse} in \cite[sec.4.4]{LZ}.

\medskip
\noindent\parbox[t]{2in}{\ul{\bf Case (ii).}\vspace{3pt} {\it The edge $e$ is an orientable loop.} \vspace{3pt}\\
In this case, the contraction \makebox(10,6){\qquad\quad increases} \vspace{4pt}}\qquad
\parbox{4in}{$G = 
  \risS{-17}{due-i2}{\put(31,8){\mbox{$e$}}\put(37,27){\mbox{$e$}}}{70}{0}{0}\quad 
   \risS{-4}{toto}{}{30}{0}{0}\quad 
  \risS{-8}{contr-i1}{}{70}{35}{15} = G/e\ .$}
 the number of vertices by 1, splitting the end-vertex of $e$ in two vertices in a way indicated by $e$.\vspace{-20pt}

\medskip
\noindent\parbox[t]{2in}{\ul{\bf Case (iii).}\vspace{3pt} {\it The edge $e$ is a non-orientable loop.} \vspace{3pt}\\
In this case the contraction}\qquad
\parbox{4in}{$G = 
  \risS{-17}{due-iii1}{\put(31,8){\mbox{$e$}}\put(37,27){\mbox{$e$}}}{70}{0}{0}\quad 
   \risS{-4}{toto}{}{30}{0}{0}\quad 
  \risS{-13}{contr-iii4}{}{70}{45}{30} = G/e\ .$}
 reverses the attachment of edge-ribbons on half of the end-vertex of $e$.

\medskip
Another notion of contraction of a loop is suggested in a recent preprint \cite{HM}. It requires a generalization of ribbon graphs whose vertices are allowed to be higher genera surfaces. For instance, their contraction of a loop leads to creation of a new vertex represented by the union of the old vertex-disc and the loop-ribbon. 

\medskip
The next lemma generalizes the contraction-deletion property of dual planar graphs mentioned above to arbitrary ribbon graphs.

{\bf Lemma.} \label{l:du-c/d} {\it 
Let $G$ be a (signed) ribbon graph, $E'\subset E(G)$ be a subset of edges of $G$, and $e\not\in E'$ be an arbitrary edge of $G$ which is not in $E'$. Then
$$(G/e)^{E'} = G^{E'\cup e}-e = G^{E'}/e\qquad\mbox{and}\qquad
(G-e)^{E'} = G^{E'\cup e}/e = G^{E'}-e\ .
$$}

The lemma obviously follows from the given definitions of contraction and deletion.
\hspace*{\fill}$\square$

\bigskip
\section{The Bollob\'as-Riordan polynomial}\label{s:br}

% To define the Bollob\'as-Riordan polynomial we need to introduce several parameters of  
% a ribbon graph $G$. 
\noindent Let\vspace{-15pt}
\begin{itemize}
\item[$\bullet$] $v(G) := |V(G)|$ denote the number of vertices of a ribbon graph $G$;
\item[$\bullet$] $e(G) := |E(G)|$ denote the number of edges of $G$;
\item[$\bullet$] $k(G)$ denote the number of connected components of $G$;
\item[$\bullet$] $r(G):=v(G)-k(G)$ be the {\it rank} of $G$;
\item[$\bullet$] $n(G):=e(G)-r(G)$ be the {\it nullity} of $G$;
\item[$\bullet$] $\bc(G)$ denote the number of connected components of the boundary of the surface of $G$.
\end{itemize} 
A {\it spanning subgraph} of a ribbon graph $G$ is %defined as 
a subgraph consisting of all the vertices of $G$ and a subset of the edges of $G$.
Let $\cF(G)$ denote the set of spanning subgraphs of $G$.
Clearly, $|\cF(G)| = 2^{e(G)}$. 
%For a signed ribbon graph we need one more parameter of a spanning subgraph. 
Let $e_{-}(F)$ be the number of negative edges in $F$.
Denote $\wo F = G - F$ the complement to $F$ in $G$, i.e.
the spanning subgraph of $G$ with exactly those (signed) edges of $G$ that do not belong to $F$.  Finally, let
$$s(F) := \frac{e_{-}(F)-e_{-}(\wo F)}{2}\ .$$

\begin{defn}\label{def:br}{\rm
The signed {\it Bollob\'as-Riordan polynomial} $R_G(x,y,z)$ is defined by
$$\fbox{$\displaystyle R_G(x,y,z)\ :=\ \sum_{F \in \cF(G)}
   x^{r(G)-r(F)+s(F)}
   y^{n(F)-s(F)}
   z^{k(F)-\bc(F)+n(F)} $}\ .
$$}
\end{defn}
In general this is a Laurent polynomial in $x^{1/2}$, $y^{1/2}$, and $z$.

The signed version of the Bollob\'as-Riordan polynomial was introduced in \cite{CP} (a version of it was also used in \cite{LM}). If all the edges are positive then it is obtained from the original Bollob\'as-Riordan polynomial
\cite{BR3} by a simple substitution $x+1$ for $x$ and $1$ for $w$. The variable $w$ in the original Bollob\'as-Riordan polynomial is responsible for orientability of the ribbon graph $F$. 
Note that the exponent $k(F)-\bc(F)+n(F)$ of the variable $z$ is equal to
$2k(F)-\chi(\widetilde{F})$, where $\chi(\widetilde{F})$ is the Euler characteristic of the surface $\widetilde{F}$ obtained by gluing a disc to each boundary component of $F$.
For orientable $F$, it is twice the genus of $F$.
In particular, for a planar ribbon graph $G$ (i.e. when the surface $G$ has genus zero)
the Bollob\'as-Riordan polynomial $R_G$ does not
depend on $z$. In this case, and if all the edges are positive, it is essentially
equal to the classical Tutte polynomial $T_{\Ga}(x,y)$ of the underlying graph $\Ga$ of $G$:
$$R_G(x-1,y-1,z) = T_{\Ga}(x,y)\ .$$
In \cite{Ka2} L.~Kauffman (see also \cite{GR})
introduced a generalization of the Tutte polynomial to signed graphs.
The previous relation holds for them as well.
Similarly, a specialization $z=1$ of the Bollob\'as-Riordan polynomial
of an arbitrary ribbon graph $G$ gives the (signed) Tutte polynomial of the underlying graph:
$$R_G(x-1,y-1,1) = T_{\Ga}(x,y)\,.$$
So one may think about the Bollob\'as-Riordan polynomial as a generalization of the Tutte polynomial to graphs embedded into a surface.

\begin{exa}\label{ex2} {\rm
Consider the ribbon graph $G$ from the example \ref{ex}(c) and shown on the left in the table below.
The other columns show eight possible spanning subgraphs $F$ and
the corresponding values of $k(F)$, $r(F)$, $n(F)$, $\bc(F)$, and $s(F)$.
$$
\begin{array}{c||c|c|c|c} \label{br-table}
\risS{8}{rg-ex32}{\put(30,35){\plinc}\put(50,30){\miinc}\put(50,-8){\miinc}}{65}{50}{10}
& \quad\risS{8}{rgBBB}{\put(42,30){\miinc}\put(42,-8){\miinc}}{58}{0}{0}
& \quad\risS{8}{rgBBA}{\put(42,-8){\miinc}}{58}{0}{0} 
& \quad\risS{13}{rgBAB}{\put(42,25){\miinc}}{58}{0}{0}  
& \quad\risS{13}{rgBAA}{}{58}{0}{0}\\ \hline
 (k,r,n,\bc,s) & (1,1,1,2,1) & (1,1,0,1,0) & (1,1,0,1,0) & (2,0,0,2,-1)\makebox(0,15){}
\\ \hline\hline
& \risS{8}{rg-ex32}{\put(30,35){\plinc}\put(50,30){\miinc}\put(50,-8){\miinc}}{65}{55}{10}
& \risS{8}{rgABA}{\put(30,35){\plinc}\put(50,-8){\miinc}}{65}{55}{10} 
& \risS{13}{rgAAB}{\put(30,30){\plinc}\put(50,25){\miinc}}{65}{55}{10}
& \risS{13}{rgAAA}{\put(30,30){\plinc}}{65}{55}{10}\\ \cline{2-5}
& (1,1,2,1,1) & (1,1,1,1,0) & (1,1,1,1,0) & (2,0,1,2,-1)\makebox(0,15){}
\end{array}
$$

We have
$$R_G(x,y,z) = x+2+y+xyz^2+2yz+y^2z\ .$$
}
\end{exa}

\begin{subs}\label{ss:prop} {\bf Properties.}\end{subs}

The Bollob\'as-Riordan polynomial is multiplicative with respect to the operations of the {\it disjoint union} $G_1\sqcup G_2$ and the {\it one-point join}
$G_1\cdot G_2$:
$$R_{G_1\sqcup G_2} = R_{G_1\cdot G_2} = R_{G_1}\cdot R_{G_2}\ .$$
Note that the operation of one-point join is ambiguous. So the equality claims that the
Bollob\'as-Riordan polynomial does not detect this ambiguity.
For unsigned ribbon graphs, these properties were proved in \cite{BR3}. For signed graphs, the proof is practically the same and follows from additivity of $s(F)$ with respect to either of these operations.

\begin{prop}\label{prop:c-d}
{\bf The signed contraction-deletion property}.\\
Let $G$ be a signed ribbon graph. Then for every positive
edge $e$ of $G$ %we have
\begin{equation}
\label{eqn:b-r_skein_pos}
\begin{array}{ll}
R_G=R_{G/e}+R_{G-e} & 
         \mbox{if $e$ is ordinary, that is neither a bridge nor a loop,} \\
R_{G}=(x+1)R_{G/e} & \mbox{if $e$ is a bridge.}
\end{array}\end{equation}
Also, for every negative edge $e$ of $G$ %we have
\begin{equation}
\label{eqn:b-r_skein_neg}
\begin{array}{ll}
R_G=x^{-1/2}y^{1/2}R_{G - e} + x^{1/2}y^{-1/2}R_{G / e} & \mbox{if $e$ is ordinary,}\\ R_{G}=x^{1/2}y^{-1/2}(y+1)R_{G/e} & \mbox{if $e$ is a bridge.}
\end{array}\end{equation}
\end{prop}

The proof of the proposition is straightforward. Spanning subgraphs of $G$ which do not contain the edge $e$ are in one-to-one correspondence with spanning subgraphs of $G-e$, while spanning subgraphs of $G$ containing $e$ are in one-to-one correspondence with spanning subgraphs of $G/e$. The equations (\ref{eqn:b-r_skein_pos}) were proved in
\cite[theorem 1]{BR3}. The equations (\ref{eqn:b-r_skein_neg}) were found by M.~Chmutov.

B.~Bollob\'as and O.~Riordan \cite{BR3} indicated a contraction-deletion property for {\it trivial loop}. A loop $e$ is called {\it trivial} if its removal and a cut of its end-vertex-disc along a chord connecting the two end segments of $e$, increase the number of connected components of the surface. In other words, there is no path from the arc $A$ to the arc $B$ in the figures of cases (ii) and (iii) in section \ref{ss:due} outside the drawn vertex. Here is an extension of the contraction-deletion properties of a trivial loop from \cite{BR3} to signed graphs. 
{\it\begin{equation}
\label{eqn:b-r_skein_loop}
\begin{array}{ll}
R_{G}=(y+1)R_{G-e} & \mbox{if $e$ is a trivial orientable (the case (ii)) 
    positive loop,} \\
R_{G}=x^{-1/2}y^{1/2}(x+1)R_{G-e} & \mbox{if $e$ is a trivial orientable   
   negative loop,} \\
R_{G}=(yz+1)R_{G-e} & \mbox{if $e$ is a trivial non-orientable 
       (the case (iii)) positive loop,} \\
R_{G}=x^{-1/2}y^{1/2}(xz+1)R_{G-e} & \mbox{if $e$ is a trivial 
   non-orientable negative loop.}
\end{array}\hspace{-10pt}\end{equation}}

S.~Huggett and I.~Moffatt gave \cite{HM} a generalization of these properties to an arbitrary (not necessary trivial) loop. However, as we already mentioned in section \ref{ss:du-c/d}, their contraction of a loop is different from ours and creates ribbon graphs with a complicated structure on vertices. Our definition of contraction does not admit a generalization of (\ref{eqn:b-r_skein_loop}) to nontrivial orientable loops. However, for arbitrary non-orientable loops we have
{\it\begin{equation}
\label{eqn:b-r_skein_no-loop}
\begin{array}{ll}
R_{G}=R_{G-e}+ yz R_{G/e} & \mbox{if $e$ is a non-orientable 
   positive loop,} \\
R_{G}=x^{-1/2}y^{1/2}(R_{G-e}+ xz R_{G/e}) & \mbox{if $e$ is a   
   non-orientable negative loop.}
\end{array}\hspace{-10pt}\end{equation}}
Note that equations (\ref{eqn:b-r_skein_no-loop}) imply the last two equations
of (\ref{eqn:b-r_skein_loop}) because if $e$ is a trivial non-orientable loop, then
both graphs $G-e$ and $G/e$ are two different one-point joins of the same two graphs. Therefore their Bollob\'as-Riordan polynomials are equal to each other.

\medskip
Our main theorem \ref{th:m-th} implies a generalization of (\ref{eqn:b-r_skein_loop}) to nontrivial orientable loops for a specialization of the Bollob\'as-Riordan polynomial to $xyz^2=1$ (see lemma \ref{le:loop} below).

\begin{prop}\label{prop:ch-of-s}
{\bf Change of the sign function}.\\
Let $G_\ve$ be a ribbon graph with the sign function $\ve$ and $G_{-\ve}$ be the same ribbon graph only with the sign function $-\ve$. Then 
$$R_{G_{-\ve}}(x,y,z) = 
    \Bigl(\frac{y}{x}\Bigr)^{(n(G_\ve)-r(G_\ve))/2}\cdot R_{G_\ve}(y,x,z)\ .
$$
\end{prop}

{\bf Proof.}
A spanning subgraph $F$ of $G_\ve$ may be regarded as a spanning subgraph of $G_{-\ve}$.
Let $s_\ve(F)$ and $s_{-\ve}(F)$ be the values of the parameter $s(\cdot)$ in graphs
$G_\ve$ and $G_{-\ve}$ respectively.
We have
$$s_{-\ve}(F) = \frac{(e(F)-e_{-}(F))-(e(\wo{F})-e_{-}(\wo{F}))}{2}
   = \frac{e(F)-e(\wo{F})}{2} - s_\ve(F)
   = e(F)-e(G_\ve)/2- s_\ve(F)\ .
$$
Hence, for the corresponding monomial of $R_{-\ve}(x,y,z)$ we get
$$\begin{array}{l}
x^{r(G_{-\ve})-r(F)+s_{-\ve}(F)} y^{n(F)-s_{-\ve}(F)} = 
x^{r(G_\ve)-r(F)+e(F)-e(G_\ve)/2-s_\ve(F)} y^{n(F)-e(F)+e(G_\ve)/2+s_\ve(F)}  
                         \vspace{8pt}\\
\qquad\qquad =x^{r(G_\ve)+n(F)-e(G_\ve)/2-s_\ve(F)} y^{-r(F)+e(G_\ve)/2+s_\ve(F)}
                         \vspace{8pt}\\
\qquad\qquad = (x^{r(G_\ve)-e(G_\ve)/2} y^{-r(G_\ve)+e(G_\ve)/2})\cdot 
            (y^{r(G_\ve)-r(F)+s_\ve(F)} x^{n(F)-s_\ve(F)}) \vspace{8pt}\\
\qquad\qquad = \Bigl(\frac{y}{x}\Bigr)^{(n(G_\ve)-r(G_\ve))/2}\cdot 
            (y^{r(G_\ve)-r(F)+s_\ve(F)} x^{n(F)-s_\ve(F)})\ .
\end{array}
$$
The monomial in the last parentheses is exactly a monomial of $R_{G_\ve}(y,x,z)$.
In other words, 
$$R_{G_{-\ve}}(x,y,z) = 
    \Bigl(\frac{y}{x}\Bigr)^{(n(G_\ve)-r(G_\ve))/2}\cdot R_{G_\ve}(y,x,z)\ .
$$
\hspace*{\fill}$\square$

\bigskip
\section{Main result}\label{s:mr}

\begin{thm}\label{th:m-th} 
The restriction of the polynomial 
$x^{k(G)} y^{v(G)} z^{v(G)+1} R_G(x,y,z)$ to the surface $xyz^2=1$ is invariant under the generalized duality. 
In other words, for any choice of the subset of edges $E'$, if $G':=G^{E'}$, then
$$x^{k(G)} y^{v(G)} z^{v(G)+1} R_G(x,y,z)\Bigl|_{xyz^2=1}\Bigr. =\ 
  x^{k(G')} y^{v(G')} z^{v(G')+1} R_{G'}(x,y,z)\Bigl|_{xyz^2=1}\Bigr. 
$$
\end{thm}
If the Bollob\'as-Riordan polynomial $R_G(x,y,z)$ contains half-integer exponents, then
the restriction to the surface $xyz^2=1$ should rather be understood as a restriction to the surface $x^{1/2}y^{1/2}z=1$.

{\bf Remark.} I.~Moffatt noticed \cite{Mof3} that in the case of orientable ribbon graphs this theorem follows from \cite[theorem 4.3]{Mof} by equating HOMFLY polynomials of some appropriate links in thickened surfaces. Thus in this case the duality relation can be derived from the link theory.

\bigskip
\begin{subs}\label{ss:proof} {\bf Proof of the theorem.}\end{subs}
By definition \ref{def:br} we may represent the polynomial
$x^{k(G)} y^{v(G)} z^{v(G)+1} R_G(x,y,z)$ as sum of monomials $M_G(F)$ corresponding to the spanning subgraphs $F$:
$$x^{k(G)} y^{v(G)} z^{v(G)+1} R_G(x,y,z)\ =\ \sum_{F \in \cF(G)} M_G(F)\ ,$$
where
$$\begin{array}{crl}
M_G(F)&:=& x^{k(G)+r(G)-r(F)+s(F)}
   y^{v(G)+n(F)-s(F)}
   z^{v(G)+1+k(F)-\bc(F)+n(F)}\\
&=&x^{k(F)+s(F)}
   y^{e(F)+k(F)-s(F)}
   z^{2k(F)+e(F)-\bc(F)+1}\\
&=&(xyz^2)^{k(F)} x^{s(F)}y^{e(F)-s(F)}z^{e(F)-\bc(F)+1}\ .
\end{array}
$$

For a spanning subgraph $F$ of $G$ we define a spanning subgraph $F'$ of $G'$ by the following rule.
\begin{quotation}\it
An edge $e$ of $G'$ belongs to the spanning subgraph $F'$ if and only if either\\
$e\in E'$ and $e\not\in F$, or $e\not\in E'$ and $e\in F$.
\end{quotation}
The correspondence $\cF(G)\ni F\leftrightarrow F'\in \cF(G')$ is one-to-one. Therefore it is enough to prove that
$$M_G(F)\Bigl|_{xyz^2=1}\Bigr. =\ M_{G'}(F')\Bigl|_{xyz^2=1}\Bigr.\ ,
$$
which is equivalent to
\begin{equation}\label{eq:mon}
x^{s(F)}y^{e(F)-s(F)}z^{e(F)-\bc(F)+1}\Bigl|_{xyz^2=1}\Bigr. =\ 
x^{s(F')}y^{e(F')-s(F')}z^{e(F')-\bc(F')+1}\Bigl|_{xyz^2=1}\Bigr.\ .
\end{equation}
By lemma \ref{le:du-prop} it is sufficient to consider the case when $E'$ consists 
of a single edge $e$, so $G'=G^{\{e\}}$.
Moreover, we may assume that $e\in F$, and hence $e\not\in F'$. Indeed, if $e\not\in F$ then $e\in F'$ and by lemma \ref{le:du-prop} $G=G'^{\{e\}}$. Therefore, interchanging
$G$ and $G'$ allows to make such an assumption.

We need to compare the parameters $s(\cdot)$,  $e(\cdot)$, and $\bc(\cdot)$ for the subgraphs $F$ and $F'$. The correspondence $F\leftrightarrow F'$ is chosen in such a way that $\bc(F)=\bc(F')$. Moreover, $e(F)=e(F')+1$, by assumption. Now, if $\ve_G(e)=+1$
then the edge $e$ does not make any contribution to $s(F)$, but in this case 
$\ve_{G'}(e)=-1$ and $e\in \wo{F'}$. Therefore $s_G(F)=s_{G'}(F')+1/2$. Similarly, if
$\ve_G(e)=-1$, then also $s_G(F)=s_{G'}(F')+1/2$. Thus we have
$$x^{s(F)}y^{e(F)-s(F)}z^{e(F)-\bc(F)+1} = x^{s(F')}y^{e(F')-s(F')}z^{e(F')-\bc(F')+1}
\cdot x^{1/2}y^{1/2}z\ ,
$$ 
which readily implies equation (\ref{eq:mon}).
\hspace*{\fill}$\square$

\begin{lemma}\label{le:loop}
{\bf Contraction-deletion of a nontrivial orientable loop}.\\
Let $e$ be an nontrivial orientable loop of a signed ribbon graph $G$. Then
\begin{equation}\label{eqn:b-r_skein_oloop}
\begin{array}{ll}
R_{G}\Bigl|_{xyz^2=1} =R_{G-e}\Bigl|_{xyz^2=1} +\ (y/x)R_{G/e}\Bigl|_{xyz^2=1} & 
       \mbox{if $e$ is positive,} \vspace{8pt}\\
R_{G}\Bigl|_{xyz^2=1} = yz\Bigl(R_{G-e} + R_{G/e}\Bigr)\Bigl|_{xyz^2=1} & 
       \mbox{if $e$ is negative.} 
\end{array}\end{equation}
\end{lemma}

{\bf Proof.}
Consider the dual graph $G^{\{e\}}$ (see the case (ii) of section \ref{ss:due}). It has one more vertex than $G$. The main theorem \ref{th:m-th} implies that 
$$R_{G}\Bigl|_{xyz^2=1} =\ yz R_{G^{\{e\}}}\Bigl|_{xyz^2=1}\ .
$$
Since $e$ was a nontrivial loop in $G$, it becomes an ordinary edge in the graph $G^{\{e\}}$. Therefore we can apply proposition \ref{prop:c-d}. If $e$ was positive in $G$, it becomes negative in $G^{\{e\}}$, and we should use the first equation of 
(\ref{eqn:b-r_skein_neg}). If $e$ was negative in $G$, it becomes positive in $G^{\{e\}}$, and we should use the first equation of (\ref{eqn:b-r_skein_pos}). In the former case, we have
$$R_{G}\Bigl|_{xyz^2=1} =\ yz \Bigl(x^{-1/2}y^{1/2}R_{G^{\{e\}}-e} +
   x^{1/2}y^{-1/2}R_{G^{\{e\}}/e}\Bigr)\Bigl|_{xyz^2=1}\ .
$$
But $G^{\{e\}}-e =G/e$ and $G^{\{e\}}/e =G-e$ according to the lemma of section \ref{ss:du-c/d} on page \pageref{l:du-c/d}. Then the substitution $x^{1/2}y^{1/2}z=1$
gives the first equation of (\ref{eqn:b-r_skein_oloop}). Similarly, the latter case when $e$ was a negative loop of $G$ implies the second equation of (\ref{eqn:b-r_skein_oloop}).
\hspace*{\fill}$\square$

\section{Natural duality of graphs on surfaces}

\begin{subs}\label{ss:BR-du} {\bf Natural duality for the Bollob\'as-Riordan polynomial.} \end{subs}

Let $G$ be a signed ribbon graph with a sign function $\ve$. The duality with respect to the set of all edges 
$E'=E(G)$ gives the dual signed graph $G_{-}^* := G^{E'}$ whose sign function is $-\ve$.
Let $G^*:=G_{+}^*$ be the graph obtained from $G_{-}^*$ by flipping the signs of all edges. In other words, $G^*$ is the same ribbon graph as $G_{-}^*$ only with the signed function changed back to $\ve$. We call $G^*$
the {\it natural dual} to $G$ because they are embedded into the same surface $\wt{G}=\wt{G^*}$ in a naturally dual way.

{\bf Proposition} (\cite[theorem 4.7]{EMS,Mof}). \ \\
{\it Let $g=k(G)-\chi(\wt{G})/2$, where $\chi(\wt{G})$ is the Euler characteristic of the surface $\wt{G}$ (it is equal to the genus of the surface in the orientable case). Then 
$$x^g R_G(x,y,z)\Bigl|_{xyz^2=1}\Bigr. = 
  y^g R_{G^*}(y,x,z)\Bigl|_{xyz^2=1}\Bigr.\ .
$$}

In \cite{EMS,Mof} this proposition was proved only for unsigned orientable ribbon 
graphs, however I.~Moffatt noticed that a combination of his theorems 4.7 and 4.3 from
\cite{Mof} implies the statement for signed orientable ribbon 
graphs as well. Here we prove it for not necessarily orientable signed graphs.
 
{\bf Proof.} The main theorem \ref{th:m-th} implies that
$$\begin{array}{rcl}
R_G(x,y,z)\Bigl|_{xyz^2=1}\Bigr. &=& 
   \left( x^{k(G_{-}^*)-k(G)} y^{v(G_{-}^*)-v(G)} z^{v(G_{-}^*)-v(G)} 
   R_{G_{-}^*}(x,y,z)\right)\Bigl|_{xyz^2=1}\Bigr.  \vspace{10pt}\\
&=& \!\! \displaystyle\Bigl(\,\frac{y}{x}\Bigr)^{(v(G^*)-v(G))/2} 
   R_{G_{-}^*}(x,y,z)\Bigl|_{xyz^2=1}\Bigr.\ ,
\end{array}
$$
because the generalized duality preserves the number of connected components,
$k(G)=k(G_{-}^*)$ (lemma \ref{le:du-prop}(f)).

The change of sign property \ref{prop:ch-of-s} gives
$$R_{G_{-}^*}(x,y,z) = 
    \Bigl(\frac{y}{x}\Bigr)^{(n(G^*)-r(G^*))/2}\cdot R_{G^*}(y,x,z)\ .
$$
Combing this with the previous equation we get
$$R_G(x,y,z)\Bigl|_{xyz^2=1}\Bigr. = \Bigl(\frac{y}{x}\Bigr)^{(v(G^*)-v(G)+n(G^*)-r(G^*))/2}\cdot 
    R_{G^*}(y,x,z)\Bigl|_{xyz^2=1}\Bigr.\ .
$$
The numerator of the exponent of this equation may be transformed to
$$\begin{array}{l}
  v(G^*)-v(G)+n(G^*)-r(G^*) =
  v(G^*)-v(G)+(e(G^*)-v(G^*)+k(G^*))-(v(G^*)-k(G^*)) \vspace{5pt}\\ 
  \hspace{3cm}
= -v(G^*) + e(G^*) - v(G) + 2k(G^*) = 2k(G) - \chi(\wt{G}) = 2g\ ,
\end{array}
$$
where $\chi(\wt{G})$ is the Euler characteristic of the surface $\wt{G}$.
\hspace*{\fill}$\square$

\begin{subs}\label{ss:Tu-du} {\bf Duality for the Tutte polynomial of planar graphs.} \end{subs}

Let $G$ be a connected plane ribbon graph, i.e. its underlying graph $\Ga$ is embedded into the plane. We assume that all edges of $G$ are positive. 
Then, as said before, the 
Bollob\'as-Riordan polynomial $R_G$ does not depend on $z$ and 
$T_{\Ga}(x,y) = R_G(x-1,y-1,z)$. The duality with respect to the set of all edges 
$E'=E(G)$ gives the usual plane dual graph $G_{-}^*=G^{E'}$ with the underlying graph $\Ga^*$ which is also connected and embedded into the same plane. 
However, all edges of $G_{-}^*$ become negative. Let
$G^*:=G_{+}^*$ be the graph obtained from $G_{-}^*$ by changing all the edges from negative to positive. Of course, $R_{G^*}(x-1,y-1,z) = T_{\Ga^*}(x,y)$. Proposition \ref{ss:BR-du} implies that
$$T_{\Ga}(x,y) = R_G(x-1,y-1,z)\Bigl|_{(x-1)(y-1)z^2=1}\Bigr. =
       R_{G^*}(y-1,x-1,z)\Bigl|_{(x-1)(y-1)z^2=1}\Bigr. = T_{\Ga^*}(y,x)\ ,
$$
because $g=0$ for the planar case.
Thus the famous duality relation for the Tutte polynomial of planar connected graphs 
$$T_{\Ga}(x,y) = T_{\Ga^*}(y,x)$$
is a direct consequence of our generalized duality for the Bollob\'as-Riordan polynomial.

\section{Virtual links and the Jones polynomial}\label{s:vl}

We will follow the Kauffman approach \cite{Ka3} to virtual links and the Jones polynomial. Such links are represented by diagrams similar to ordinary link diagrams, except some crossings are designated as {\it virtual}. Virtual crossings should be understood not as
crossings but rather as defects of our two-dimensional figures. They should
be treated in the same way as the extra crossings appearing in planar pictures of
non-planar graphs. Here are some examples of virtual knots.
$$\risS{-18}{ex}{}{65}{20}{20}\hspace{3cm}
  \risS{-18}{v31}{}{40}{0}{20}\hspace{3cm}
  \risS{-18}{v41}{}{55}{0}{20}\label{ex:vir-kn}
$$
On figures we encircle the virtual crossings to distinguish them from the classical ones.

Virtual link diagrams are considered up to plane isotopy, the {\it classical}
Reidemeister moves:% involving classical crossings:
$$\risS{-10}{RI}{}{65}{20}{12}\qquad\qquad
  \risS{-10}{RII}{}{65}{0}{0}\qquad\qquad
  \risS{-10}{RIII}{}{65}{0}{0}\quad ,
$$
and the {\it virtual} Reidemeister moves:
$$\risS{-10}{RI-v}{}{65}{17}{15}\qquad\quad
  \risS{-10}{RII-v}{}{65}{0}{0}\qquad\quad
  \risS{-10}{RIII-v}{}{65}{0}{0}\qquad\quad
  \risS{-10}{RIV-v}{}{65}{0}{0}\quad .
$$

We define the Jones polynomial for virtual links using the Kauffman bracket in the same way as for classical links. Let $L$ be a virtual link diagram.
Consider two ways of resolving a classical crossing.
The {\it $A$-splitting},\ $\smf{cr}\ \leadsto\ \smf{Asp}$\ ,
is obtained by joining the two vertical angles swept out by the overcrossing arc when
it is rotated counterclockwise toward the undercrossing arc.
Similarly, the {\it $B$-splitting},\ $\smf{cr}\ \leadsto\ \smf{Bsp}$\ , is
obtained by joining the other two vertical angles. A {\it state} $s$ of
a link diagram $L$
is a choice of either an $A$ or $B$-splitting at each classical crossing.
Denote by $\cS(L)$ the set of states of $L$.
A diagram $L$ with $n$ crossings has $|\cS(L)| = 2^n$
different states.

Denote by $\a(s)$ and $\b(s)$ the numbers of $A$-splittings and $B$-splittings
in a state $s$, respectively, and by $\d(s)$ the number of
components of the curve obtained from the link
diagram $L$ by 
splitting according to the state $s \in \cS(L)$. Note that virtual crossings do not connect components.

\begin{defn}\label{def:kb}
The \emph{Kauffman bracket} of a diagram $L$ is a polynomial in three variables
$A$, $B$, $d$ defined by the formula
$$ \kb{L} (A,B,d)\ :=\ \sum_{s \in \cS(L)} \,
A^{\a(s)} \, B^{\b(s)} \, d^{\,\d(s)-1}\,.
$$
\end{defn}

Note that $\kb{L}$ is \emph{not} a topological
invariant of the link; it depends on the link
diagram. However, it determines the \emph{Jones polynomial}
$J_L(t)$ by a simple substitution
$$A=t^{-1/4},\qquad B=t^{1/4},\qquad d=-t^{1/2}-t^{-1/2}\ ;$$
$$J_L(t)\, := (-1)^{w(L)} t^{3w(L)/4} \kb{L} (t^{-1/4}, t^{1/4}, -t^{1/2}-t^{-1/2})\ .
$$
Here $w(L)$ denotes the {\it writhe} determined
by an orientation of $L$ as the sum over the classical
crossings of $L$ of the signs\,: 
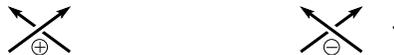
\vspace{-10pt}
$$\risS{-10}{cr_pl}{}{25}{0}{15}\hspace{3cm}\label{loc-writhe}
  \risS{-10}{cr_mi}{}{25}{0}{0}\quad .
$$
The Jones polynomial is a topological invariant (see e.g.~\cite{Ka1}).

For example, for virtual knots above the Kauffman bracket and the Jones polynomial are
the following.
$$\begin{array}{c@{\qquad}l}
\risS{-22}{ex}{}{65}{0}{0} & \kb{L} = A^3 + 3A^2Bd + 2AB^2 + AB^2d^2+ B^3d \vspace{6pt}\\
    & J_L(t)=1 \vspace{15pt}\\
\risS{-25}{v31}{}{40}{0}{0} &  \kb{L} = A^2d + 2 AB + B^2\vspace{6pt}\\
    & J_L(t) = t^{-3/2} + t^{-1} - t^{-1/2} \vspace{25pt}\\
\risS{-25}{v41}{}{55}{0}{0} & \kb{L} = A^3d + 3A^2B + 2AB^2 + AB^2d+ B^3d \vspace{6pt}\\
    & J_L(t) = t^{-2}-t^{-1}-t^{-1/2}+1+t^{1/2}
\end{array}\vspace{25pt}$$

\section{Thistlethwaite's type theorems}

In 1987 Thistlethwaite \cite{Th} (see also \cite{Ka1}) proved that up to a sign and a power of $t$, the Jones polynomial 
$V_L(t)$ of an alternating link $L$ is equal to 
the Tutte polynomial $T_{\Ga_L}(-t,-t^{-1})$
of the planar graph $\Ga_l$ obtained from a checkerboard coloring of the
regions of a link diagram. 
$$\risS{0}{kd-gr}{\put(5,50){\mbox{$L$}}  \put(240,45){\mbox{$\Ga_L$}}
                 }{250}{70}{10} 
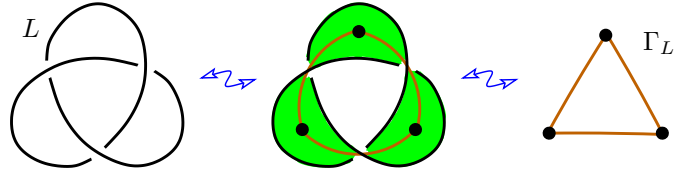
\label{p:thist}
$$

L.~Kauffman \cite{Ka2} generalized the theorem to arbitrary (classical) links using signed graphs.

\begin{subs}\label{ss:virt-Th} {\bf Thistlethwaite's theorem for virtual links.} \end{subs}

Here we explain a generalization of this theorem to virtual
links. With each state $s$ of a virtual link diagram $L$ we associate a ribbon graph
(possibly non-orientable) $G_L^s$. 
We express the Kauffman bracket (and hence the Jones polynomial) of $L$ as a specialization of the Bollob\'as-Riordan polynomial of $G_L^s$. For two different states
$s$ and $s'$, the ribbon graphs $G_L^s$ and $G_L^{s'}$ are dual to each other. Then
the generalized duality for the Bollob\'as-Riordan polynomial implies that the result of the specialization does not depend on the choice of state $s$. 

Our construction of $G_L^s$ is a straightforward generalization of the construction from \cite{DFKLS} where it was used for classical links only.
The vertices of $G_L^s$ are obtained by gluing a disc to each state circle of $s$. Let us describe the edges of $G_L^s$. In a vicinity of a classical crossing of $L$ we place a small \vyd{planar} band connecting two strands of the splitting of $s$. These bands will be the edge-ribbons of $G_L^s$. The orientation of the plane induces the orientations on these bands and since the arrows on the state circles.
These arrows indicate how the edge-ribbons are attached to the vertices. By the definition \ref{def:rb2} this information determines the ribbon graph $G_L^s$.
If a crossing of $L$ is resolved as an $A$-splitting in the state $s$, we assign $+1$ to the corresponding edge, if it is resolved as 
a $B$-splitting, then we assign $-1$. So we get a sign function. 

The next example illustrates this construction. Here we enumerate the crossings and assign the same numbers to the bands and arrows corresponding to them.
$$\begin{array}{l}
\risS{-20}{constr-chvo1}{\put(-3,40){\mbox{$L$}}
  \pn{21}{31}{$1$}\pn{40}{15}{$2$}\pn{60}{16}{$3$}\pn{15}{-10}{Diagram}
  \pn{125}{32}{$1$}\pn{145}{14}{$2$}\pn{166}{15}{$3$}\put(115,43){\splinc}
      \put(145.5,32){\smiinc}\put(175,21){\smiinc}\pn{125}{-10}{State $s$}
  \pn{233}{31}{$1$}\pn{254}{13}{$2$}\pn{276}{15}{$3$}\put(224,44){\splinc}
      \put(253.5,31){\smiinc}\put(283,21){\smiinc}
      \pn{210}{-10}{Attaching planar bands}
  \pn{333}{40}{$1$}\pn{347}{48}{$1$}\pn{366}{13}{$2$}\pn{372}{33}{$2$}
      \pn{390}{15}{$3$}\pn{400}{25}{$3$}\pn{315}{-10}{Replacing bands by arrows}
  }{400}{40}{90} \\
\hspace{1cm}\risS{-20}{constr-chvo2}{
         \put(130,60){\mbox{1 ---\ \ \splinc\ \ ; 2 ---\ \ \smiinc\ \ ; 3 ---\ \ \smiinc}}
  \pn{33}{31}{$1$}\pn{48}{48}{$1$}\pn{74}{13}{$2$}\pn{80}{34}{$2$}
      \pn{100}{15}{$3$}\pn{111}{25}{$3$}\pn{24}{-10}{Untwisting state circles}
  \pn{175}{9}{$1$}\pn{175}{45}{$1$}\pn{158}{28}{$2$}\pn{224}{44}{$2$}
      \pn{194}{28}{$3$}\pn{224}{8}{$3$}\pn{150}{-7}{Pulling state circles apart}
  \pn{288}{55}{$1$}\pn{348}{30}{$2$}\pn{348}{14}{$3$}
      \put(324,55){\splinc}\put(355,47){\smiinc}\put(355,-2){\smiinc}
      \pn{275}{-12}{Forming the ribbon graph $G_L^s$}
  }{380}{0}{45} 
\end{array}$$

\begin{lemma}\label{le:two-st}
Let $s$ and $s'$ be two states of the same diagram $L$. Then the graphs $G_L^s$ and 
$G_L^{s'}$ are dual with respect to a set of edges corresponding to the crossings where the states $s$ and $s'$ are different from each other.
\end{lemma}

The proof of the lemma is straightforward. We exemplify it by a figure of the construction of the ribbon graph $G_L^{s'}$.

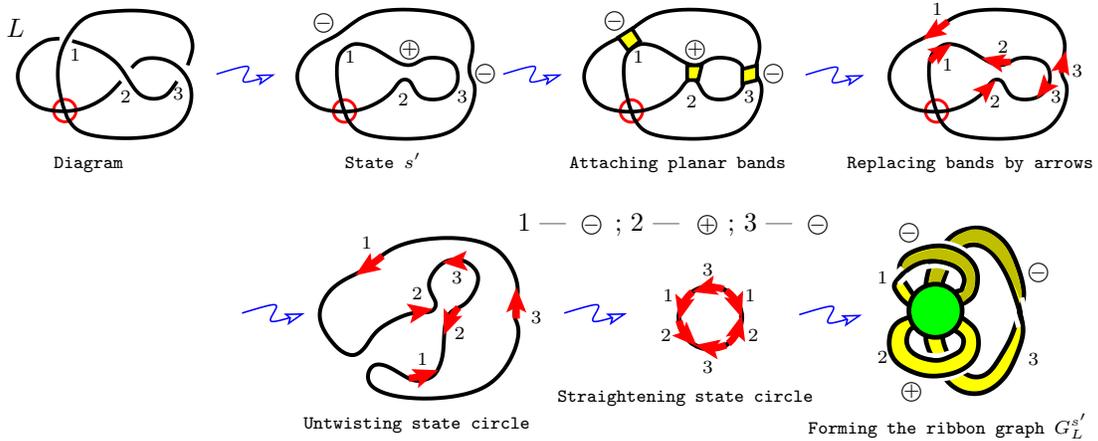
$$\begin{array}{l}
\risS{-20}{constr-rgst-pr}{\put(-3,40){\mbox{$L$}}
  \pn{21}{31}{$1$}\pn{40}{15}{$2$}\pn{60}{16}{$3$}\pn{15}{-10}{Diagram}
  \pn{125}{30}{$1$}\pn{145}{14}{$2$}\pn{167.5}{15}{$3$}\put(114,43){\smiinc}
      \put(145.5,33){\splinc}\put(175,24){\smiinc}\pn{125}{-10}{State $s'$}
  \pn{234}{30}{$1$}\pn{254}{13}{$2$}\pn{276}{15}{$3$}\put(224,44){\smiinc}
      \put(253.5,33){\splinc}\put(284,24){\smiinc}
      \pn{210}{-10}{Attaching planar bands}
  \pn{351}{27}{$1$}\pn{347}{48}{$1$}\pn{369}{13}{$2$}\pn{372}{34}{$2$}
      \pn{392}{15}{$3$}\pn{400}{25}{$3$}\pn{315}{-10}{Replacing bands by arrows}
  }{400}{40}{90} \\
\hspace{3cm}\risS{-20}{constr-rgst-pr2}{
         \put(105,65){\mbox{1 ---\ \ \smiinc\ \ ; 2 ---\ \ \splinc\ \ ; 3 ---\ \ \smiinc}}
  \pn{67}{14}{$1$}\pn{46}{58}{$1$}\pn{65}{39}{$2$}\pn{81}{24}{$2$}
      \pn{80}{45}{$3$}\pn{110}{30}{$3$}\pn{24}{-10}{Untwisting state circle}
  \pn{160}{38}{$1$}\pn{191}{38}{$1$}\pn{159}{23}{$2$}\pn{192}{23}{$2$}
      \pn{175}{11}{$3$}\pn{175}{48}{$3$}\pn{120}{0}{Straightening state circle}
  \pn{241}{45}{$1$}\pn{241}{15}{$2$}\pn{298}{14}{$3$}
      \put(250,2){\splinc}\put(299,47){\smiinc}\put(250,62){\smiinc}
      \pn{215}{-12}{Forming the ribbon graph $G_L^{s'}$}
  }{300}{0}{45} 
\end{array}$$
In this example $G_L^{s'}=(G_L^s)^{\{1,2\}}$ (see page \pageref{p:ex-c}).
\hspace*{\fill}$\square$

\begin{thm}\label{th:v-th}
Let $L$ be a virtual link diagram with $e$ classical crossings, $G_L^s$ be the signed ribbon graph corresponding to a state $s$, and
$v:=v(G_L^s)$, $k:=k(G_L^s)$. Then $e=e(G_L^s)$ and
$$\kb{L} (A,B,d) = A^e\left( x^k y^v z^{v+1} 
    R_{G_L^s}(x,y,z)\Bigl|_{x=\frac{Ad}{B},\ y=\frac{Bd}{A},\ z=\frac{1}{d}}\Bigr. 
    \right)\ .\vspace{5pt}
$$
\end{thm}
Note, that the substitution $x=\frac{Ad}{B}$, $y=\frac{Bd}{A}$, and $z=\frac{1}{d}$
satisfies the equation $xyz^2=1$. Then lemma \ref{le:two-st} and the main theorem
\ref{th:m-th} (the duality property of the Bollob\'as-Riordan polynomial) imply that
the right-hand side of the equation does not depend on the initial state $s$. Therefore it is enough to proof the theorem for one particular choice of $s$. Let us impose an orientation on the diagram $L$ and pick a state $s$, {\it Seifert state}, where all the splittings respect orientations of strands. Then our construction of the graph $G_L^s$
literally coincides with the construction from \cite{CV}. The sign function in this case 
is equal to the local writhe of a crossing from p.\pageref{loc-writhe}. 
Hence the following theorem from
\cite{CV} implies our theorem.
\hspace*{\fill}$\square$

\begin{thm}[{\cite[theorem 4.1]{CV}}] \label{th:cv}
Let $L$ be a virtual link diagram, $s$ be the Seifert state, $G_L^s$ be the corresponding signed ribbon graph, and
$n:=n(G_L^s)$, $r:=r(G_L^s)$, $k:=k(G_L^s)$. Then
$$\kb{L} (A,B,d) = A^n B^r d^{k-1} \,
R_{G_L^s}\left(\frac{Ad}{B}, \frac{Bd}{A}, \frac{1}{d}\right)\ .
$$
\end{thm}

It is also not difficult to prove theorem \ref{th:v-th} directly by making a one-to one correspondence between states of $L$ and spanning subgraphs of $G_L^s$ and proving the equality of the corresponding monomials. In this correspondence, the initial state $s$ should correspond to the spanning subgraph of $G_L^s$ without any edges. The dual state
$\hat{s}$, that is the one with all opposite splittings as compared to $s$, corresponds to the whole graph $G_L^s$. In general, the number of boundary components of a spanning subgraph of $G_L^s$ should be equal to the number of state circles of the corresponding state.

\begin{cor}\label{cor:cv} 
For a state $s$ of an oriented virtual link diagram $L$ with $e$ classical crossings we have
$$J_L(t) = (-1)^{w(L)} t^{(3w(L)-e+2r)/4} (-t^{1/2}-t^{-1/2})^{k-1}
R_{G_L^s}\left(-1-t^{-1}, -t-1, \frac{1}{-t^{1/2}-t^{-1/2}}\right)\ .
$$
\end{cor}

\bigskip
\begin{subs}\label{ss:constr-dfkls} {\bf The theorem from \cite{DFKLS}.}\end{subs}

In the paper \cite{DFKLS}, the authors dealt with classical links only. They used the state
$s_A$ consisting of $A$-splittings only. In this case, their construction of the ribbon graph gives our  $G_L^{s_A}$. It is always orientable for classical links. All the edges are positive, so we do not need the signed versions of the polynomials here. Also they used a slightly different versions of the Kauffman bracket and the Bollob\'as-Riordan polynomial related to ours as follows:
$$\langle L \rangle(A) := [L](A,A^{-1},-A^2-A^{-2})\ ,\qquad\qquad
  C(G;\ X,Y,Z) := R_G(X-1,Y,Z^{1/2})\ .
$$

\begin{thm}[{\cite[theorem 5.4]{DFKLS}}] \label{th:dfkls}
For a connected classical link diagram $L$, 
$$\langle L \rangle(A) = A^{e+2-2v} \,
C\left(G_L^{s_A};\ -A^4, -1-A^{-4}, \frac{1}{(-A^2-A^{-2})^2}\right)\ ,
$$
where $e=e(G_L^{s_A})$is the number of crossings of $L$, and $v=v(G_L^{s_A})$ is the number of vertices of the Turaev surface $G_L^{s_A}$.
\end{thm}

The connection of $L$ implies the connection of the ribbon graph, $k(G_L^{s_A})=1$. Therefore, this theorem is a special case of theorem \ref{th:v-th} when $s=s_A$ and $L$ is classical and connected.
\hspace*{\fill}$\square$

Gluing a disk to each boundary component of $G_L^s$ we obtained a closed surface without boundary, $\wt{G_L^s}$, which is called the {\it Turaev surface} of the state $s$.
V.~Turaev used the surface $\wt{G_L^{s_A}}$ in his proof of the Tait conjectures \cite{Tu}. The projection of the link $L$ onto the Turaev surface is alternating. 
The same is true for virtual links, however the notion of an alternating link diagram on the
non-orientable Turaev surface requires a clarification.

\bigskip
\begin{subs}\label{ss:constr-cp} {\bf The theorem from \cite{CP}.}\end{subs}

The paper \cite{CP} is devoted to another version of Thistlethwaite's theorem for a particular class of virtual link diagrams known as {\it checkerboard colorable} diagrams.
This notion was introduced by N.~Kamada in \cite{Kam1, Kam2}, who showed that many classical results on knots and links can be extended to checkerboard colorable virtual links. A similar notion under the name of {\it atom} was studied by V.~Manturov  \cite{Man} following A.~Fomenko \cite{Fom}.

Checkerboard colorability is related to a {\it coorientation}.
A {\it coorientation} (see, for example, \cite{Ar}) of a plane curve is a choice of one of the two sides of the curve in a neighborhood of each point of the curve. We depict the chosen side by indicating a normal direction field along the curve. 
A link diagram $L$ is called  
{\it checkerboard colorable} if there is a coorientation of a state $s$ of $L$ such that
near each crossing the coorientations of the two strands point to the opposite directions.
$$\risS{-10}{coor-cr1}{}{35}{0}{0}\qquad\qquad\risS{-10}{coor-cr2}{
      \pn{-95}{-10}{Possible coorientations near a crossing}}{35}{30}{0}
\hspace{4cm}
\risS{-10}{coor-cr3}{}{35}{0}{0}\qquad\qquad\risS{-10}{coor-cr4}{
      \pn{-110}{-10}{Change of coorientations of strands at a crossing}}{35}{0}{25}
$$
It is easy to see that if this condition is satisfied for a state $s$ then it will be satisfied for any other state as well. Also, we may think about coorientation of the original diagram $L$ which is changing to the opposite one when the point passes a crossing (see the right two figures above). In other words, the coorientation of strands of $L$ near a crossing has to point into two vertical angles. Of course, near a virtual crossing the coorientation of one strand goes through without noticing the crossing strand or its coorientation.

All classical link diagrams are checkerboard colorable since we may take the normal direction field pointing inside the green regions (see the trefoil figure on page \pageref{p:thist}) of the checkerboard coloring of the regions of the diagram.
Among the three virtual knot diagrams on page \pageref{ex:vir-kn}, only the first one is 
checkerboard colorable.
$$\risS{-18}{ex-coor}{
      \pn{0}{-10}{Checkerboard colorable}}{75}{35}{40}\hspace{3cm}
  \risS{-15}{v31-41-coor}{
      \pn{50}{-6}{Error}
      \pn{40}{-20}{Not checkerboard colorable}}{180}{0}{0}
$$
It is easy to see that all alternating virtual link diagrams are checkerboard colorable.
For such diagrams there is a {\it canonical} checkerboard coloring \cite{Kam2} when near each crossing the coorienting field points inside vertical angles that are glued together by the $A$-splitting. In the figure above, the right most figure has a fragment of a canonical coloring near the crossing.
N.~Kamada proved \cite{Kam1} that a virtual link diagram $L$ is checkerboard
colorable if and only if it can be made alternating by switching some classical overcrossings to undercrossings. Let $L^{alt}$ be the alternating diagram obtained from $L$ in this way. For $L^{alt}$, we consider the canonical checkerboard coloring. It determines a state $s_B$ by performing $B$-splittings at all classical crossings of $L^{alt}$. The unsigned ribbon graph $G_L^{CP}$ constructed in 
\cite[sections 3 and 4]{CP} is the same as $G_{L^{alt}}^{s_B}$ constructed in section \ref{ss:virt-Th}. However, the sign function $\ve$ for $G_L^{CP}$ is different.
If an edge $e$ of $G_L^{CP}$ corresponds to a crossing where the switching was performed during the way from $L$ to $L^{alt}$, then we set $\ve(e)=-1$, and we set $\ve(e)=+1$ for the other edges. 
\begin{thm}[{\cite[theorem 4.1]{CP}}] \label{th:cp}
For a checkerboard colorable virtual link diagram $L$, 
$$\kb{L} (A,B,d)  \ = \ A^{r(G_L^{CP})} B^{n(G_L^{CP})} d^{k(G_L^{CP})-1} \,
R_{G_L^{CP}}\left(\frac{Bd}{A}, \frac{Ad}{B}, \frac{1}{d}\right).
$$
\end{thm}

{\bf Proof.}
Let $s$ be a state of $L$ given by $A$-splitting of the crossings that were switched 
during the way from $L$ to $L^{alt}$ and by $B$-splitting of the crossings that were
not switched. The ribbon graph $G_{L^{alt}}^{s_B}=G_L^{CP}$ is the same as $G_L^s$. Moreover the sign function $\ve$ for $G_L^{CP}$ is precisely opposite to the signed function for $G_L^s$ defined in the section \ref{ss:virt-Th}. By the change of sign property
\ref{prop:ch-of-s} we have
$$R_{G_L^{CP}}(x,y,z) = 
    \Bigl(\frac{y}{x}\Bigr)^{(n(G_L^s)-r(G_L^s))/2}\cdot R_{G_L^s}(y,x,z)\ .
$$
Let $v:=v(G_L^s)$, $k:=k(G_L^s)$, $e:=e(G_L^s)$, $r:=r(G_L^s)=v-k$, and
$n:=n(G_L^s)=e-v+k$.

By theorem \ref{th:v-th} we have
$$\kb{L} (A,B,d) = A^e 
  \Bigl(\frac{Ad}{B}\Bigr)^k \Bigl(\frac{Bd}{A}\Bigr)^v \Bigl(\frac{1}{d}\Bigr)^{v+1} 
    R_{G_L^s}\left(\frac{Ad}{B},\frac{Bd}{A},\frac{1}{d}\right)\ .
$$
Combing the last two equations together we get
$$\begin{array}{rcl}
  \kb{L} (A,B,d) &=& \displaystyle A^e 
  \Bigl(\frac{Ad}{B}\Bigr)^k \Bigl(\frac{Bd}{A}\Bigr)^v \Bigl(\frac{1}{d}\Bigr)^{v+1}
  \Bigl(\frac{Bd/A}{Ad/B}\Bigr)^{(n-r)/2}
    R_{G_L^{CP}}\left(\frac{Bd}{A}, \frac{Ad}{B}, \frac{1}{d}\right) \vspace{10pt}\\
&=&\displaystyle A^{e+k-v-n+r} B^{-k+v+n-r} d^{k-1}
    R_{G_L^{CP}}\left(\frac{Bd}{A}, \frac{Ad}{B}, \frac{1}{d}\right)  \vspace{10pt}\\
&=&\displaystyle A^r B^n d^{k-1}
    R_{G_L^{CP}}\left(\frac{Bd}{A}, \frac{Ad}{B}, \frac{1}{d}\right)\ .
\end{array}
$$
\hspace*{\fill}$\square$

So the theorem of \cite{CP} follows from the theorem of \cite{CV} by using the generalized duality theorem \ref{th:m-th}. 

\begin{rem}\rm
Theorems of \cite{CP} and \cite{DFKLS} were also unified in a recent preprint \cite{Mof2} whose construction of the {\it unsigning} (see section 3.1 there) is a special case of our generalized duality.
\end{rem}

\section{Possible further directions}

{\bf 1.} B.~Gr\"unbaum and G.~Shephard \cite{GS} define duality as a bijection from the sets of vertices and faces of one polyhedron to the sets of faces and vertices of another polyhedron preserving incidence. They noticed that the square of the self-duality map might not be equal to the identity map. Further, D.~Archdeacon and R.~Richter \cite{AR} described all spherical self-dual polyhedra. The same definition may be applied to general ribbon graphs. For example, the graph on a torus from the Introduction (the same as in example \ref{ex}) is self-dual with respect to the set of  all edges. Also a ribbon graph consisting of a single vertex and a single non-orientable loop is self-dual as a graph embedded into the projective plane. It would be interesting to investigate the duality
maps for ribbon graphs with respect to a set of edges. Do the constructions of \cite{AR} give all self-dual ribbon graphs for higher genera as well? 

\bigskip
{\bf 2.} There are various parameters \cite{A,Thom} which measure a density of the embedding of a ribbon graph $G$ into a closed surface $\wt{G}$ without boundary. For instance, the {\it edge-width}, $ew(G)$, is the length of the shortest noncontractible cycle in $G$; the {\it face-width}, $fw(G)$, is the minimum of the number of points of intersection $C\cap\Ga$ of a noncontractible cycle $C$ in $\wt{G}$ with the underlying graph $\Ga$ of $G$ taken over all such $C$; the {\it dual-width}, $dw(G)$, is the minimum of the number of points of intersection $C\cap\Ga$ taken over noncontractible cycles $C$ in $\wt{G}$ which intersect $\Ga$ only in the interior of its edges. It would be interesting to explore the behavior of these parameters with respect to the generalized duality and also their relations to the Bollob\'as-Riordan polynomial. 

\bigskip
{\bf 3.} J.~Edmonds \cite{Ed} found a condition on two abstract graphs with a bijection between the sets of their edges to be embedded into the same surface in a naturally dual way. The condition is formulated in terms of the bijection. The problem is to find a similar condition for the generalized duality with respect to a subset of edges.

\bigskip
{\bf 4.} The Bollob\'as-Riordan polynomial of $G$, as a generalization of the Tutte polynomial of its underlying abstract graph $\Ga$, shares with the latter its various remarkable specializations. For instance, see \cite{B} for details, for a connected ribbon graph $G$ we have
\begin{itemize}
\item $T_{\Ga}(1,1)=R_G(0,0,1)$ is the number of spanning trees in $G$;
\item $T_{\Ga}(2,1)=R_G(1,0,1)$ is the number of edge sets forming forests in $G$;
\item $T_{\Ga}(1,2)=R_G(0,1,1)$ is the number of connected spanning subgraphs in $G$;
\item $T_{\Ga}(2,2)=R_G(2,1,1)$ is the number of spanning subgraphs in $G$.
\end{itemize}
Also, for any ribbon graph $G$ the {\it chromatic polynomial} $\chi_{\Ga}(\lambda)$ of its underlying graph $\Ga$ equals
$$\chi_{\Ga}(\lambda)= (-1)^{r(G)}\lambda^{k(G)}T_{\Ga}(1-\lambda,0) =
(-1)^{r(G)}\lambda^{k(G)}R_G(-\lambda,-1,1)\ .
$$
M.~Korn and I.~Pak \cite{KP} found a combinatorial meaning of the specialization
$R_G(k,k,1/k)$ for any natural number $k$, which is expressed in terms of the numbers of monochromatic vertices of and edge-coloring of $\Ga$  in $k$ colors subject to certain restrictions. The $k=2$ case is related to the number of T-{\it tetromino tilings} of the
torus ribbon graphs.
Note that this specialization satisfies the equation $xyz^2=1$. According to our main theorem \ref{th:m-th} such specializations of $G$ and its generalized dual $G^{E'}$ are proportional to each other. It would be interesting to find a direct combinatorial
bijective proof as well as other combinatorial interpretations of various specializations of the Bollob\'as-Riordan polynomial.

\bigskip
{\bf 5.} The Bollob\'as-Riordan polynomial has a multivariable version 
\cite[section 2.2]{HM} where each edge is labeled by its own variable. This version generalizes the multivariable Tutte polynomial from \cite{S}. It would be interesting to extend the generalized duality to such labeled graphs and find a multivariable analog of our main theorem \ref{th:m-th}.
\footnote{After the paper was submitted to the journal this was done by Fabien Vignes-Tourneret \cite{VT}.}

\vskip1.cm

\parbox[t]{2.5in}{\it \textbf{Sergei~Chmutov}\\
Department of Mathematics\\
The Ohio State University, Mansfield\\
1680 University Drive\\
Mansfield, OH 44906\\
~\texttt{chmutov@math.ohio-state.edu}} 

\end{document}